\documentclass[12pt]{amsart}
\newtheorem{thm}{Theorem}[section]
\newtheorem{cor}[thm]{Corollary}
\newtheorem{lem}[thm]{Lemma}
\newtheorem{prop}[thm]{Proposition}
\theoremstyle{definition}

\theoremstyle{remark}

\numberwithin{equation}{section}
\begin{document}

\begin{center}
{\bf{Curvature properties of some class of warped product manifolds}}
\end{center}

\vspace{1mm}

\begin{center}
Ryszard Deszcz, Ma\l gorzata G\l ogowska, Jan Je\l owicki and Georges Zafindratafa
\end{center}

\vspace{1mm}

\begin{center}
{\sl{Dedicated to the memory of Professor W\l odzimierz Waliszewski}}
\end{center}

\vspace{1mm}

\noindent
{\bf{Abstract.}}
{\footnotesize{Warped product manifolds with p-dimensional base, p=1,2,
satisfy some curvature conditions of pseudosymmetry type.
These conditions are formed from the metric tensor g, the Riemann-Christoffel curvature tensor R, 
the Ricci tensor S and the Weyl conformal curvature C of the considered manifolds. 
The main result of the paper states that 
if p=2 and the fibre is a semi-Riemannian space of constant curvature, if n is greater or equal to 4,
then the (0,6)-tensors
R.R - Q(S,R) and C.C of such warped products
are proportional to the (0,6)-tensor Q(g,C) 
and the tensor C is expressed by a linear combination of some Kulkarni-Nomizu products 
formed from the tensors g and S. Thus these curvature conditions satisfy non-conformally flat non-Einstein
warped product spacetimes (p=2, n=4).
We also investigate curvature properties of pseudosymmetry type of quasi-Einstein manifolds.   
In particular, we obtain some curvature property of the Goedel spacetime.}}\footnote{{\bf{Mathematics 
Subject Classification (2010):}} 
Primary 53B20, 53B25, 53B30, 53B50, 83C15; Secondary 53C25, 53C40, 53C80.

{\bf{Key words and phrases:}}
Einstein manifold, quasi-Einstein manifold, three-dimensional Berger sphere, 
three-spheres of Kaluza-Klein type,
2-quasi-Einstein manifold, warped product manifold,
spacetime, Robertson-Walker spacetime, generalized Robertson-Walker spacetime,
G\"{o}del spacetime, general spherically symmetric spacetime, Vaidya spacetime,
hypersurface, 2-quasi-umbilical hypersurface, Chen ideal submanifold,
pseudosymmetric manifold, manifold with pseudosymmetric Weyl tensor, 
peudosymmetry type curvature condition.

The first named author is supported by the
Universit\'{e} de Valenciennes et du Hainaut-Cambr\'{e}sis, France
and by a grant of the Wroc\l aw University of Environmental and Life Sciences,
Poland.
The second and third named authors are supported
by a grant of the Wroc\l aw University of Environmental and Life Sciences, Poland.}

\section{Introduction}

%\hspace*{6mm}
Let $g$, $\nabla$, $R$, $S$, $\kappa $ and $C$ be the metric tensor, the Levi-Civita connection, 
the Riemann-Christoffel curvature tensor, the Ricci tensor, the scalar curvature tensor 
and the Weyl conformal curvature tensor 
of a semi-Riemannian manifold $(M,g)$, $n = \dim M \geq 3$, respectively. 
It is well-known that $(M,g)$ is said 
to be an {\sl Einstein manifold} if at every point of $M$ 
its Ricci tensor $S$ is proportional to the metric tensor $g$, 
i.e., $S = \frac{\kappa}{n}\, g$ on $M$ \cite{Besse}.
In particular, if $S = 0$ on $M$ then $(M,g)$ is called a {\sl Ricci flat manifold}. 
We denote by ${\mathcal U}_{S}$ the set of all points of $(M,g)$ at which 
$S$ is not proportional to $g$, i.e.,
${\mathcal U}_{S} \, = \,  \{x \in M\, | \, 
S - \frac{\kappa }{n}\, g \neq 0\ \mbox {at}\ x \}$.
The manifold $(M,g)$ is said to be 
a {\sl quasi-Einstein manifold} if 
\begin{eqnarray}
\mathrm{rank}\, (S - \alpha\, g) &=& 1
\label{quasi02}
\end{eqnarray}
on ${\mathcal U}_{S}$, where $\alpha $ is some function on ${\mathcal U}_{S}$.
In particular, if $\mathrm{rank}\, S = 1$ on ${\mathcal U}_{S}$ then  
$(M,g)$ is called a {\sl Ricci-simple manifold} \cite{DRV}. 
Every warped product manifold $\overline{M} \times _{F} \widetilde{N}$
of an $1$-dimensional $(\overline{M}, \overline{g})$ base manifold and
a $2$-dimensional manifold $(\widetilde{N}, \widetilde{g})$
or an $(n-1)$-dimensional Einstein manifold
$(\widetilde{N}, \widetilde{g})$, $n \geq 4$, with a warping function $F$,
is a  quasi-Einstein manifold (see, e.g., {\cite[Section 1] {ChDGP}}). 
We mention that quasi-Einstein manifolds arose during the study of exact solutions
of the Einstein field equations and the investigation on quasi-umbilical hypersurfaces 
of conformally flat spaces, see, e.g., \cite{DGHSaw} and references therein. 
Quasi-Einstein hypersurfaces in semi-Riemannian spaces of constant curvature
were studied among others in:
\cite{{DeGl}, {DGHS}, {R102}, {DHS105}, {G6}}, see also \cite{DGHSaw}.
We refer to \cite{{ChDGP}, {DeHoJJKunSh}} 
for recent results on quasi-Einstein manifolds.
The semi-Riemannian manifold $(M,g)$, $n \geq 3$, 
is called a {\sl $2$-quasi-Einstein manifold} if 
\begin{eqnarray}
\mathrm{rank}\, (S - \alpha \, g ) &\leq & 2 ,
\label{quasi0202weak}
\end{eqnarray}
on ${\mathcal U}_{S}$ and $\mathrm{rank}\, (S - \alpha \, g ) = 2$
on some open non-empty subset of ${\mathcal U}_{S}$, 
where $\alpha $ is some function on ${\mathcal U}_{S}$ (see, e.g., \cite{{DGHZ02}, {DGP-TV02}}).
Every warped product manifold $\overline{M} \times _{F} \widetilde{N}$
of a $2$-dimensional base manifold $(\overline{M}, \overline{g})$ 
and a $2$-dimensional manifold $(\widetilde{N}, \widetilde{g})$
or an $(n-2)$-dimensional Einstein semi-Riemannian manifold
$(\widetilde{N}, \widetilde{g})$, $n \geq 5$, 
with a warping function $F$, satisfies (\ref{quasi0202weak})
(see Theorem 6.1 of this paper).
Some exact solutions of the Einstein field equations are 
non-conformally flat $2$-quasi-Einstein manifolds.
For instance, the Reissner-Nordstr\o m spacetime, as well as
the Reissner-Nordstr\o m-de Sitter type spacetimes are such manifolds (see, e.g.,
\cite{Kow02}). It seems that the Reissner-Nordstr\o m spacetime
is the "oldest" example of 
a non-conformally flat 
$2$-quasi-Einstein warped product manifold.
It is easy to see that every $2$-quasi-umbilical hypersurface
in a semi-Riemannian space of constant curvature is a $2$-quasi-Einstein manifold
(see, e.g., \cite{DGP-TV02}). 

Let $A$ and $B$ be symmetric $(0,2)$-tensors on a semi-Riemannian manifold $(M,g)$.
We denote by $A \wedge B$ their Kulkarni-Nomizu tensor. 
We note that (\ref{quasi02}) 
holds at a point $x \in {\mathcal U}_{S} \subset M$ if and only if at this point we have
$(S - \alpha\, g) \wedge (S - \alpha\, g) = 0$, i.e.  
\begin{eqnarray}
\frac{1}{2}\, S \wedge S   - \alpha\, g \wedge S + \alpha ^{2}\, G &=& 0, \ \ G \ =\ \frac{1}{2}\, g \wedge g .    
\label{quasi03}
\end{eqnarray}
From (\ref{quasi03}), by a suitable contraction,
we get immediately 
\begin{eqnarray}
S^{2} &=& (\kappa - (n-2) \alpha )\, S + \alpha ( (n-1) \alpha - \kappa )\, g.
\label{quasi03quasi03} 
\end{eqnarray}
For precise definitions of the symbols used here, 
we refer to Section 2 of this paper (see also \cite{{ChDGP}, {DGHSaw}}). 
We can write the Weyl conformal curvature tensor $C$ of $(M,g)$, $n \geq 3$, by 
\begin{eqnarray}
C &=& R - \frac{1}{n-2}\, g\wedge S + \frac{\kappa }{(n-2)(n-1)}\, G .
\label{eqn2.1}   
\end{eqnarray}
It is well-known that a semi-Riemannian manifold $(M,g)$, $n \geq 4$, is conformally flat 
if and only if $C = 0$ everywhere in $M$. From $C = 0$, by (\ref{eqn2.1}), we get immediately
\begin{eqnarray}
R &=&  \frac{1}{n-2}\, g \wedge S - \frac{\kappa }{(n-2) (n-1)} \, G . 
\label{dddWeyl}
\end{eqnarray}
The Robertson-Walker spacetimes, and more generally, warped products of an $1$-dimen\-sio\-nal manifold 
and an $(n-1)$-dimensional semi-Riemannian space of constant curvature, $n \geq 4$,
are conformally flat quasi-Einstein manifolds (see, e.g., {\cite[Section 4] {K 3}}). 
It is obvious that 
(\ref{quasi03}) and (\ref{dddWeyl}) yield 
\begin{eqnarray*}
R &=& \frac{1}{2}\, S \wedge S   + \left( \frac{1}{n-2} - \alpha \right) g \wedge S 
+ \left( \alpha ^{2} - \frac{\kappa }{(n-2) (n-1)} \right) G 
%\label{quasi04}
\end{eqnarray*}
(see, e.g., {\cite[p. 150] {DGHHY}}). 
Thus the curvature tensor $R$ of a conformally flat quasi-Einstein manifold $(M,g)$, $n \geq 4$,
is expressed by a linear combination of the tensors: $S \wedge S$, $g \wedge S$ and $G$. 
We also can investigate non-conformally flat and non-quasi-Einstein semi-Riemannian manifolds $(M,g)$, $n \geq 4$, 
whose curvature tensor $R$ is a linear combination
of these tensors.
More precisely, we can investigate semi-Riemannian manifolds $(M,g)$, $n \geq 4$,
satisfying on the set ${\mathcal U}_{S} \cap {\mathcal U}_{C} \subset M$ the condition
\begin{eqnarray}
R &=& \frac{\phi}{2}\, S\wedge S + \mu\, g\wedge S + \eta\, G ,
\label{eq:h7a}
\end{eqnarray}
where 
${\mathcal U}_{C}$ is the set of all points of $M$ at which $C \neq 0$ and
$\phi$, $\mu $ and $\eta $ are some functions on ${\mathcal U}_{S} \cap {\mathcal U}_{C}$.
A semi-Riemannian manifold $(M,g)$, $n \geq 4$, satisfying (\ref{eq:h7a}) on 
${\mathcal U}_{S} \cap {\mathcal U}_{C} \subset M$ 
is called a {\sl Roter type manifold, or Roter manifold, or Roter space} 
\cite{{P106}, {DGP-TV01}, {DGP-TV02}}. 
Roter type manifolds and in particular Roter type hypersurfaces 
in semi-Riemannian spaces of constant curvature were studied in:
\cite{{P106}, {DGHHY}, {DGHZ01}, {DGHZ02}, {DGP-TV01}, {R102}, {DeKow}, {DePlaScher}, {DeScher}, {G108}, {G5}, {Kow01}, {Kow02}}. 
In Section 3 we present curvature conditions satisfying by Roter type manifolds.
In particular, on every Roter type manifold $(M,g)$, $n \geq 4$, 
the following relations are satisfied  
on ${\mathcal U}_{S} \cap {\mathcal U}_{C} \subset M$: 
\begin{eqnarray}
R \cdot R - Q(S,R) &=& L\, Q(g,C) ,
\label{genpseudo01}\\          
C \cdot C &=& L_{C}\, Q(g,C) ,
\label{4.3.012}\\
C \cdot R + R \cdot C 
&=& Q(S,C) + \left( L + L_{C} - \frac{1}{ (n-2) \phi } \right) Q(g,C) ,
\label{newRoter02}\\ 
C \cdot R - R \cdot C 
&=& Q(S,C) - \frac{\kappa }{n-1}\, Q(g,C) ,
\label{newRoter01}
\end{eqnarray}
where 
$L = L_{R} + \frac{\mu }{\phi }$,
$L_{C} = L_{R} - \frac{ \kappa }{n-1} + \frac{1}{(n-2)\phi} - \frac{\mu}{\phi}$ and
$L_{R} = \frac{1}{\phi} ( (n-2) ( \mu ^{2} - \phi \eta ) - \mu )$
(Theorem 3.2 and Proposition 3.3).
In 
{\cite[Theorem 3.2] {S3}}
(also see {\cite[Section 4] {DGP-TV02}} and {\cite[Section 4] {Saw06}})
it was proved that the curvature tensor $R$ of some hypersurfaces 
in semi-Riemannian spaces of constant curvature 
is a linear combination of the tensors:
$S \wedge S$, $g \wedge S$, $G$, $g \wedge S^{2}$, $S \wedge S^{2}$ and $S^{2} \wedge S^{2}$.
Precisely, 
we have on ${\mathcal U}_{S} \cap {\mathcal U}_{C} \subset M$ 
\begin{eqnarray}
R &=& \frac{\phi _{1}}{2}\, S \wedge S + \phi _{2}\, g \wedge S + \phi _{3}\, G + \phi _{4}\, g \wedge S^{2}
+ \phi _{5}\, S \wedge S^{2} + \frac{ \phi _{6}}{2}\, S^{2} \wedge S^{2} , 
\label{eqn0101.01}
\end{eqnarray}
where $\phi _{1}, \phi _{2}, \ldots , \phi _{6}$
are some functions on this set. 
Evidently, (\ref{eq:h7a}) is a special case of (\ref{eqn0101.01}).
Examples of manifolds satisfying (\ref{eqn0101.01}) are given in 
{\cite[Example 2.1] {P140}}, {\cite[Section 4] {DGP-TV02}},
{\cite[Example 4.1] {DeHoJJKunSh}}, {\cite[Section 5] {Saw06}} and {\cite[Section 5] {SDHJK}}.
Manifolds satisfying (\ref{eqn0101.01}) were studied in \cite{{DGJP-TZ01}, {G5}, {ShaKun01}, {ShaKun02}}.

It is easy to verify that on any semi-Riemannian manifold $(M,g)$, $n \geq 4$, the following identity is satisfied
\begin{eqnarray}
C \cdot R  + R \cdot C 
&=& R \cdot R + C \cdot C - \frac{1}{(n-2)^{2}}\, Q(g,  - \frac{\kappa }{ n-1}\, g\wedge S + g \wedge S^{2} ) 
\label{identity01}
\end{eqnarray}
(Theorem 3.4(i)). In addition, 
if (\ref{genpseudo01}), with some function $L$, holds on ${\mathcal U}_{C} \subset M$ then (\ref{identity01}) turns into 
\begin{eqnarray}
C \cdot R  + R \cdot C 
&=& Q(S,C) + L\, Q(g,C) + C \cdot C \nonumber\\
& & - \frac{1}{(n-2)^{2}}\, Q(g, \frac{n-2}{2}\, S \wedge S - \kappa \, g\wedge S + g \wedge S^{2} ) 
\label{02identity01}
\end{eqnarray}
(Theorem 3.4(ii)). 
Moreover, if (\ref{4.3.012}), with some functions $L_{C}$, is satisfied on ${\mathcal U}_{C} \subset M$
then (\ref{02identity01}) takes the form 
\begin{eqnarray}
C \cdot R + R \cdot C &=& Q(S,C) + ( L + L_{C} )\, Q(g,C)\nonumber\\
& & - \frac{1}{(n-2)^{2}}\, Q( g, 
\frac{n-2}{2}\, S \wedge S - \kappa \, g \wedge S + g \wedge S^{2} ) 
\label{identity05}
\end{eqnarray}
(Theorem 3.4(iii)). 
We note that if $(M,g)$ is a quasi-Einstein semi-Riemannian manifold 
satisfying (\ref{quasi02}) then  (\ref{identity05}),
by making use of (\ref{quasi03}) and (\ref{quasi03quasi03}), turns into
\begin{eqnarray}
C \cdot R + R \cdot C  &=& Q(S,C) + ( L + L_{C} )\, Q(g,C) ,
\label{identity05quasi}
\end{eqnarray}
and in particular, if $(M,g)$ is the G\"{o}del spacetime then (\ref{identity05quasi}) yields 
\begin{eqnarray}
C \cdot R + R \cdot C &=& Q(S,C) + \frac{\kappa }{6}\, Q(g,C) 
\label{Goedelidentity05}
\end{eqnarray}
(Theorem 3.4(iv)-(v)). The conditions (\ref{genpseudo01}) and (\ref{4.3.012})
are also satisfied on some submanifolds isometrically immersed in an Euclidean space, as well as on some hypersurfaces 
isometrically immersed in a semi-Riemannian space of constant curvature (theorems 3.7-3.9).

In Section 4 we prove that warped product manifolds  
$\overline{M} \times _{F} \widetilde{N}$ 
of an $1$-dimensional semi-Riemannian manifold $(\overline{M},\overline{g})$ and
some $(n-1)$-dimensional semi-Riemannian manifold $(\widetilde{N},\widetilde{g})$, $n \geq 4$, 
satisfy 
(\ref{genpseudo01}), (\ref{4.3.012}) and (\ref{identity05}) 
(theorems 4.1-4.3).
In particular, we state that the warped product of an $1$-dimensional manifold $(\overline{M},\overline{g})$ and
some $3$-dimensional Riemannian manifold:
the $3$-dimensional Berger spheres, 
the Heisenberg group $Nil_{3}$, 
$\widetilde{PSL(2,{\mathbb{R}})}$ - the universal covering of the Lie group $PSL(2,{\mathbb{R}})$, 
the Lie group $Sol_{3}$,
a Riemannian manifold isometric to an open part of the $3$-dimensional Cartan hypersurface
or some three-spheres of Kaluza-Klein type,
satisfies 
(\ref{genpseudo01}), (\ref{4.3.012}) and (\ref{identity05}) (Theorem 4.2).

In the next section we present results on pseudosymmetric warped product manifolds.
In particular, we consider warped products $\overline{M} \times _{F} \widetilde{N}$ 
of a $2$-dimensional semi-Riemannian manifold $(\overline{M},\overline{g})$ and
an $(n-2)$-dimensional semi-Riemannian manifold $(\widetilde{N},\widetilde{g})$, $n \geq 4$, 
with the warping function $F$,  
assuming that $(\widetilde{N},\widetilde{g})$ is a semi-Riemannian space of constant curvature, 
when $n \geq 5$. In Theorem 5.3 we present necessary and sufficient condition for 
such manifold to be pseudosymmetric. 

In Section 6 
we consider warped products $\overline{M} \times _{F} \widetilde{N}$ 
of a $2$-dimensional semi-Riemannian manifold $(\overline{M},\overline{g})$ and
an $(n-2)$-dimensional semi-Riemannian manifold $(\widetilde{N},\widetilde{g})$, $n \geq 4$, 
with the warping function $F$,  
assuming that $(\widetilde{N},\widetilde{g})$ is an Einstein semi-Riemannian manifold, when $n \geq 5$. 
Theorem 6.2 states that on some subset 
${\mathcal U}_{S} \cap {\mathcal U}_{C} \subset \overline{M} \times _{F} \widetilde{N}$
(see to that section for details) the tensor $R \cdot S$ is a linear combination
of the Tachibana tensors $Q(g,S)$, $Q(g,S^{2})$ and $Q(S,S^{2})$, i.e.
\begin{eqnarray}
R \cdot S &=&
\psi_{5} \, Q(g,S) + \psi_{4}\, Q(g, S^{2}) + \psi_{3}\, Q(S,S^{2}) ,
\label{RdotS01}
\end{eqnarray}
on this set, for some functions $\psi_{3}, \psi_{4}$ and $\psi_{5}$.
We mention that recently in \cite{DGP-TV02} it was shown that the tensor $R \cdot S$ 
of some minimal hypersurfaces in Euclidean spaces has this property
(see also \cite{{DGPSS}, {SawCM}}).
The condition (\ref{RdotS01}), by (\ref{identity21ddhh}), turns into
\begin{eqnarray}
C \cdot S &=&
\psi_{1} \, Q(g,S) + \psi_{2}\, Q(g, S^{2}) + \psi_{3}\, Q(S,S^{2}) ,
\label{RdotS02}
\end{eqnarray}
where $\psi_{1} = \psi_{5} +  \frac{\kappa }{(n-2)(n-1)}$ and 
$\psi_{2} =  \psi_{4} - \frac{1}{n-2}$.
Semi-Riemannian manifolds, and in particular, hypersurfaces in semi-Riemannian spaces of constant curvature,
satisfying the special cases of (\ref{RdotS02}), i.e. 
$C \cdot S = \psi \, Q(g,S)$,
resp., $C \cdot S = 0$,
were investigated, among others, 
in \cite{{DG90}, {Kow01}, {Kow02}}, 
resp., \cite{{DGHSaw}, {DGHS}, {DGHZ01}, {DGHZ02}, {DHS4}, {DHS01}, {SDHJK}}. 

In the last section we consider
warped products $\overline{M} \times _{F} \widetilde{N}$ 
of a $2$-dimensional semi-Riemannian manifold $(\overline{M},\overline{g})$ and
an $(n-2)$-dimensional semi-Riemannian manifold $(\widetilde{N},\widetilde{g})$, $n \geq 4$, 
with the warping function $F$,  
assuming that $(\widetilde{N},\widetilde{g})$ is a semi-Riemannian space of constant curvature, 
when $n \geq 5$. 
In Theorem 7.1(i) we state that
(\ref{genpseudo01}), (\ref{4.3.012}) and (\ref{identity05}) 
hold on ${\mathcal U}_{C} \subset \overline{M} \times _{F} \widetilde{N}$.
In Theorem 7.1(ii),
under some additional assumption, we state that on some open subset    
$V \subset {\mathcal U}_{S} \cap {\mathcal U}_{C} \subset \overline{M} \times _{F} \widetilde{N}$
the Weyl tensor $C$ of the considered warped product 
is a linear combination of the Kulkarni-Nomizu tensors
$S \wedge S$, $g \wedge S$, $g \wedge S^{2}$ and $G$.
Precisely, (\ref{mainTT032main}) holds on $V$. Evidently, (\ref{mainTT032main}) by (\ref{eqn2.1}) turns into (\ref{eqn0101.01}). 
Thus we have a new family of manifolds satisfying (\ref{eqn0101.01}).
On the set $({\mathcal U}_{S} \cap {\mathcal U}_{C}) \setminus V$ 
the Weyl tensor $C$ 
is a linear combination of the Kulkarni-Nomizu tensors
$S \wedge S$, $g \wedge S$ and $G$.
In that section we also present curvature properties of the Vaidya spacetime,
as well as of some generalized Vaidya spacetimes:
the Vaidya-Kottler, the Vaidya-Reissner-Nordstr\o m and the Vaidya-Bonnor spacetime.

\section{Preliminary results}

Throughout this paper all manifolds are assumed
to be connected paracompact manifolds of class $C^{\infty }$.
Let $(M,g)$ be an $n$-dimensional, $n \geq 2$,
semi-Riemannian manifold
and let $\nabla$ be its
Levi-Civita connection
and $\Xi (M)$
the Lie algebra of vector fields on $M$.
We define on $M$ the endomorphisms
$X \wedge _{A} Y$ and
${\mathcal{R}}(X,Y)$
of $\Xi (M)$, respectively, by
\begin{eqnarray*}
(X \wedge _{A} Y)Z \ =\
A(Y,Z)X - A(X,Z)Y,\ \
{\mathcal R}(X,Y)Z
\ =\
\nabla _X \nabla _Y Z - \nabla _Y \nabla _X Z - \nabla _{[X,Y]}Z ,
\end{eqnarray*}
where $A$ is a symmetric $(0,2)$-tensor on $M$
and $X, Y, Z \in \Xi (M) $.
The Ricci tensor $S$,
the Ricci operator ${\mathcal{S}}$,
the tensors $S^{2}$ and $S^{3}$ and
the scalar curvature $\kappa $
of $(M,g)$
are defined by
$S(X,Y)\, =\, \mathrm{tr} \{ Z \rightarrow {\mathcal{R}}(Z,X)Y \}$,
$g({\mathcal{S}}X,Y)\, =\, S(X,Y)$,
$S^{2}(X,Y) \, =\, S({\mathcal{S}}X,Y)$,
$S^{3}(X,Y) \, =\, S^{2}({\mathcal{S}}X,Y)$
and
$\kappa \, =\, \mathrm{tr}\, {\mathcal{S}}$, respectively.
The endomorphism
${\mathcal{C}}(X,Y)$ is defined by
\begin{eqnarray*}
{\mathcal{C}}(X,Y)Z  &=& {\mathcal{R}}(X,Y)Z
- \frac{1}{n-2}(X \wedge _{g} {\mathcal{S}}Y + {\mathcal{S}}X \wedge _{g} Y
- \frac{\kappa}{n-1}X \wedge _{g} Y)Z .
\end{eqnarray*}
Now the $(0,4)$-tensor $G$,
the Riemann-Christoffel curvature tensor $R$ and
the Weyl conformal curvature tensor $C$ of $(M,g)$ are defined by
$G(X_1,X_2,X_3,X_4)  = g((X_1 \wedge _{g} X_2)X_3,X_4)$ and
\begin{eqnarray*}
R(X_1,X_2,X_3,X_4)  \ =\ g({\mathcal{R}}(X_1,X_2)X_3,X_4) ,\ \
C(X_1,X_2,X_3,X_4)  \ =\ g({\mathcal{C}}(X_1,X_2)X_3,X_4) ,
\end{eqnarray*}
respectively, where $X_1,X_2,X_3,X_4 \in \Xi (M)$.
Let ${\mathcal{B}}$ be a tensor field sending any $X, Y \in \Xi (M)$
to a skew-symmetric endomorphism ${\mathcal{B}}(X,Y)$,
and let $B$ be a $(0,4)$-tensor
associated with ${\mathcal{B}}$ by
\begin{eqnarray}
B(X_1,X_2,X_3,X_4) &=&
g({\mathcal{B}}(X_1,X_2)X_3,X_4)\, .
\label{DS5}
\end{eqnarray}
The tensor $B$ is said to be a
{\sl{generalized curvature tensor}}
if
the following conditions are satisfied
\begin{eqnarray*}
& &
B(X_1,X_2,X_3,X_4) \ =\  B(X_3,X_4,X_1,X_2)\, ,\\
& &
B(X_1,X_2,X_3,X_4)
+ B(X_3,X_1,X_2,X_4)
+ B(X_2,X_3,X_1,X_4) \ =\  0\, .
\end{eqnarray*}
For ${\mathcal{B}}$ as above, let $B$ be again defined by (\ref{DS5}).
We extend the endomorphism
${\mathcal{B}}(X,Y)$ to a derivation
${\mathcal{B}}(X,Y) \cdot \, $
of the algebra of tensor fields on $M$,
assuming that it commutes with contractions and
$\ {\mathcal{B}}(X,Y) \cdot \, f \, =\, 0$, for any smooth function $f$ on $M$.
For a $(0,k)$-tensor field $T$, $k \geq 1$,
we can define the $(0,k+2)$-tensor $B \cdot T$ by
\begin{eqnarray*}
& & (B \cdot T)(X_1,\ldots ,X_k,X,Y) \ =\
({\mathcal{B}}(X,Y) \cdot T)(X_1,\ldots ,X_k)\\
&=& - T({\mathcal{B}}(X,Y)X_1,X_2,\ldots ,X_k)
- \cdots - T(X_1,\ldots ,X_{k-1},{\mathcal{B}}(X,Y)X_k)\, .
\end{eqnarray*}
In addition, if $A$ is a symmetric $(0,2)$-tensor then we define
the $(0,k+2)$-tensor $Q(A,T)$ by
\begin{eqnarray*}
& & Q(A,T)(X_1, \ldots , X_k, X,Y) \ =\
(X \wedge _{A} Y \cdot T)(X_1,\ldots ,X_k)\\
&=&- T((X \wedge _A Y)X_1,X_2,\ldots ,X_k)
- \cdots - T(X_1,\ldots ,X_{k-1},(X \wedge _A Y)X_k)\, .
\end{eqnarray*}
The tensor $Q(A,T)$
is called the {\sl Tachibana tensor of the tensors} $A$ and $T$,
or shortly the Tachibana tensor 
(see, e.g., \cite{{DGHHY}, {DGHZ01}, {P140}, {DGP-TV02}, {DeHoJJKunSh}}).
For a symmetric $(0,2)$-tensor $E$ and a $(0,k)$-tensor $T$, $k \geq 2$,
we define their Kulkarni-Nomizu product $E \wedge T$ by (\cite{DG90})
\begin{eqnarray*}
& &(E \wedge T )(X_{1}, \ldots , X_{4}, Y_{3}, \ldots , Y_{k})\\
&=&
E(X_{1},X_{4}) T(X_{2},X_{3}, Y_{3}, \ldots , Y_{k})
+ E(X_{2},X_{3}) T(X_{1},X_{4}, Y_{3}, \ldots , Y_{k} )\\
& &
- E(X_{1},X_{3}) T(X_{2},X_{4}, Y_{3}, \ldots , Y_{k})
- E(X_{2},X_{4}) T(X_{1},X_{3}, Y_{3}, \ldots , Y_{k}) .
\end{eqnarray*}
For instance,
the following tensors are generalized curvature tensors:
$R$, $C$, $G$ and $E \wedge F$,
where
$E$ and $F$ are symmetric $(0,2)$-tensors.
For a symmetric $(0,2)$-tensor $A$ we define
the endomorphism ${\mathcal{A}}$
and the tensors $A^{2}$ and  $A^{3}$
by
$g({\mathcal{A}}X,Y) = A(X,Y)$,
$A^{2}(X,Y) = A( {\mathcal{A}}X, Y)$
and
$A^{3}(X,Y) = A^{2}( {\mathcal{A}}X, Y)$,
respectively.
Let
$B_{hijk}$, $T_{hijk}$, and $A_{ij}$ be the local
components of
generalized curvature tensors
$B$ and $T$ and
a symmetric $(0,2)$-tensor
$A$ on $M$, respectively, where
$h,i,j,k,l,m,p,q \in \{ 1,2, \ldots , n \}$.
The local components
$(B \cdot T)_{hijklm}$ and
$Q(A,T)_{hijklm}$ of the tensors
$B \cdot T$, $Q(A,T)$, $B \cdot A$ and $Q(g,A)$
are the following
\begin{eqnarray}
(B \cdot T)_{hijklm}
&=&
g^{pq}(
T_{pijk}B_{ qhlm}
+ T_{hpjk}B_{ qilm}
+ T_{hipk}B_{ qjlm}
+ T_{hijp}B_{ qklm}) ,\label{genformulas01}\\
Q(A,T)_{hijklm}
&=&
A_{hl}T_{ mijk}
+ A_{il}T_{ hmjk}
+ A_{jl}T_{ himk}
+ A_{kl}T_{ hijm}\nonumber\\
& &-
A_{hm}T_{ lijk}
- A_{im}T_{ hljk}
- A_{jm}T_{ hilk}
- A_{km}T_{ hijl} ,\label{genformulas02}\\
(B \cdot A)_{hklm}
&=&
g^{pq}( A_{pk} B_{qhlm} + A_{ph} B_{qklm}),\label{genformulas03}\\
Q(g,A)_{hklm}
&=&
g_{hl}A_{km} + g_{kl}A_{hm} - g_{hm}A_{kl} - g_{km}A_{hl}  .
\label{genformulas04}
\end{eqnarray}
\begin{lem} 
Let $(M,g)$, $n \geq 3$, be a semi-Riemannian manifold.
Let $A$ be a symmetric $(0,2)$-tensor on $M$ such that 
$\mathrm{rank}(A) = 2$ at some point $x \in M$.
(i) cf. {\cite[Lemma 2.1] {P106}} The tensors $A$, $A^{2}$ and $A^{3}$ satisfy at $x$ the following relations
\begin{eqnarray}
A^{3} &=& \mathrm{tr}(A)\, A^{2} + \frac{1}{2} ( \mathrm{tr}(A^{2}) - (\mathrm{tr}(A))^{2})\, A ,
\label{eqn10.1}\\ 
A \wedge A^{2} &=& \frac{1}{2} \mathrm{tr}(A)\, A \wedge A ,
\label{eqn14.1}\\
A^{2} \wedge A^{2} &=& - \frac{1}{2} ( \mathrm{tr}(A^{2}) - (\mathrm{tr}(A))^{2})\, A \wedge A ,
\label{eqn14.2}\\ 
(A^{2} - \mathrm{tr}(A)\, A) \wedge (A^{2} - \mathrm{tr}(A)\, A )
&=& - \frac{1}{2} ( \mathrm{tr}(A^{2}) - (\mathrm{tr}(A))^{2})\, A \wedge A .
\label{eqn14.2dd}
\end{eqnarray}
(ii) Let $T$ be a generalized curvature tensor on $M$ satisfying 
\begin{eqnarray}
T &=& \frac{\phi _{0}}{2}\, A \wedge A + \phi _{2}\, g \wedge A + \phi _{3}\, G + \phi _{4}\, g \wedge A^{2}
+ \phi _{5}\, A \wedge A^{2} + \frac{ \phi _{6}}{2}\, A^{2} \wedge A^{2} , 
\label{eqn14.3}
\end{eqnarray}
where $\phi _{0}$, $\phi _{2},\ \ldots \ , \phi _{6}$
are some functions on $M$. Then at given point $x$ we have  
\begin{eqnarray*}
T &=& \frac{\phi _{1}}{2}\, A \wedge A + \phi _{2}\, g \wedge A + \phi _{3}\, G + \phi _{4}\, g \wedge A^{2} ,\\
\phi _{1} &=& \phi _{0} +  \mathrm{tr}(A)\, \phi _{5}
- \frac{1}{2} ( \mathrm{tr}(A^{2}) - (\mathrm{tr}(A))^{2})\, \phi _{6} .
\end{eqnarray*}
\end{lem}
{\bf{Proof.}} (i) (\ref{eqn10.1}) and (\ref{eqn14.1}) were already obtained in  {\cite[eqs. (2.6) and (2.10)] {P106}}.
Further, transvecting equation (2.10) of \cite{P106}, i.e.
\begin{eqnarray*}
\mathrm{tr}(A)\, ( A_{il} A_{jm} - A_{im}A_{jl}) 
+  A_{jl} A^{2}_{im} + A_{im}A^{2}_{jl} - A_{il} A^{2}_{jm} - A_{jm}A^{2}_{il} &=& 0 ,
\end{eqnarray*}
with $A^{m}_{k} = g^{ms}A_{sk}$ we obtain
\begin{eqnarray*}
A^{2}_{il} A^{2}_{jk} - A^{2}_{ik}A^{2}_{jl} +  A_{il} A^{3}_{jk} - A_{jl}A^{3}_{ik}
&=& 
\mathrm{tr}(A)\, ( A_{il} A^{2}_{jk} - A_{jl}A^{2}_{ik} ),
\end{eqnarray*}
where $g_{hk}$, $g^{hk}$, $A_{hk}$, $A^{2}_{hk}$ and $A^{3}_{hk}$ 
are the local components of the tensors $g$, $g^{-1}$, $A$, $A^{2}$ and $A^{3}$,
respectively. This, by (\ref{eqn14.1}), turns into
\begin{eqnarray*}
A^{2}_{il} A^{2}_{jk} - A^{2}_{ik}A^{2}_{jl} 
&=& 
- \frac{1}{2} ( \mathrm{tr}(A^{2}) - (\mathrm{tr}(A))^{2})\,
( A_{il} A_{jk} - A_{ik}A_{jl} ),
\end{eqnarray*}
i.e. (\ref{eqn14.2}).
Now, using (\ref{eqn14.1}) and (\ref{eqn14.2}) we get immediately (\ref{eqn14.2dd}),
which completes the proof of (i). (ii) is an obvious consequence of (i). 
\begin{lem} 
Let $B$ be a symmetric $(0,2)$-tensor on a $2$-dimensio\-nal semi-Rieman\-nian
manifold $(M, g)$. 
(i) {\cite[Lemma 2(iii)] {30}} The following identity is satisfied on $M$
\begin{eqnarray}
g \wedge B &=& \mathrm{tr}(B)\, G . 
\label{eqn10.2}
\end{eqnarray}
(ii)
The following identities are satisfied on $M$
\begin{eqnarray}
B^{2} &=& \mathrm{tr}(B)\, B + \frac{1}{2} ( \mathrm{tr}(B^{2}) - (\mathrm{tr}(B))^{2})\, g ,
\label{eqn10.3}\\ 
Q(B,B^{2}) &=& - \frac{1}{2} ( \mathrm{tr}(B^{2}) - (\mathrm{tr}(B))^{2})\, Q(g,B)\nonumber . 
%\label{eqn10.4} 
\end{eqnarray}
\end{lem}
{\bf{Proof.}} (ii) From (\ref{eqn10.2}) we get
\begin{eqnarray}
g_{hk}B_{ij} +  g_{ij}B_{hk} - g_{hj}B_{ik} -  g_{ik}B_{hj} 
&=& \mathrm{tr}(B)\, ( g_{hk}g_{ij} - g_{hj}g_{ik} ) ,
\label{01new01}  
\end{eqnarray}
where $B_{ij}$ and  $B^{2}_{ij}$ 
are the local components of the tensors $B$ and $B^{2}$, respectively.
Transvecting (\ref{01new01}) with $B^{hk} = B_{ij}g^{hi}g^{kj}$ we obtain
\begin{eqnarray*}
B^{2}_{ij} &=& \mathrm{tr}(B)\, B_{ij} + \frac{1}{2} ( \mathrm{tr}(B^{2}) - (\mathrm{tr}(B))^{2})\, g_{ij} ,  
\end{eqnarray*}
i.e. (\ref{eqn10.3}). Further, we also have
\begin{eqnarray*}
Q(B,B^{2}) 
&=& 
Q(B, \mathrm{tr}(B)\, B + \frac{1}{2} ( \mathrm{tr}(B^{2}) - (\mathrm{tr}(B))^{2})\, g)
\ =\ 
- \frac{1}{2} ( \mathrm{tr}(B^{2}) - (\mathrm{tr}(B))^{2}) \, Q(g,B),  
\end{eqnarray*}
which completes the proof.
\newline

For symmetric $(0,2)$-tensors $E$ and $F$ we have 
\begin{eqnarray}
Q(E, E \wedge F) &=& - \frac{1}{2}\, Q(F, E \wedge E ) , \ \ \ \ %\label{DS7}
E \wedge Q(E,F) \ =\ - \frac{1}{2}\, Q(F, E \wedge E )
\label{DS7}
\end{eqnarray}
(see, e.g., {\cite[Section 3] {DGHS}} and {\cite[eq. (3)] {DG90}}).
In particular, from (\ref{DS7}) we obtain
\begin{eqnarray}
Q(S, g \wedge S) &=& - \frac{1}{2}\, Q(g, S \wedge S ), 
\ \ \ \
Q(g, g \wedge S) \ =\ - Q(S, G ). 
\label{DS7new}
\end{eqnarray}
Using now (\ref{eqn2.1}) and (\ref{DS7new}) we get
\begin{eqnarray}
Q(S,R) &=& Q(S,C) - \frac{1}{n-2}\, Q(g, \frac{1}{2}\, S \wedge S) - \frac{\kappa }{ (n-2)(n-1) }\, Q(S,G) . 
\label{identity02}
\end{eqnarray}
We also have 
\begin{eqnarray}
& & 
( g\wedge S) \cdot ( g\wedge S) \ =\ - Q( S^{2} , G) ,\ \ 
G \cdot ( g\wedge S) \ =\  Q( g , g \wedge S) \ = \ - Q( S, G) ,
\label{identity21}\\
& &
( g\wedge S) \cdot S \ =\  Q( g, S^{2}) ,\ \ 
G \cdot S \ =\  Q( g,S) 
\label{identity21dd}
\end{eqnarray}
(see, e.g., {\cite[Lemma 2.1 (ii)] {DGHHY}} and {\cite[Lemma 3.2] {Kow02}}).
Using (\ref{eqn2.1}) and (\ref{identity21dd}) we obtain 
\begin{eqnarray}
C \cdot S &=& R \cdot S - \frac{1}{n-2} \, Q(g, S^{2}) +  \frac{\kappa }{(n-2)(n-1)}\, Q(g, S)   
\label{identity21ddhh}
\end{eqnarray}
(see, e.g., {\cite[p. 217] {DHS01}}).

\section{Some curvature conditions}

A semi-Riemannian manifold $(M,g)$, $n \geq 3$, is called {\sl semisymmetric} 
if $R \cdot R = 0$ on $M$ \cite{Sz 1}.
A semi-Riemannian manifold $(M,g)$, $n \geq 3$, is said to be {\sl pseudosymmetric} 
if the tensors $R \cdot R$ and $Q(g,R)$ are linearly dependent at every point of $M$
\cite{{30}, {D 10}, {38}, {DG1987}}. 
This is equivalent 
on ${\mathcal{U}}_{R} = \{x \in M\, | \, R - \frac{\kappa }{(n-1)n}\, G \neq 0\ \mbox {at}\ x \}$ to 
\begin{eqnarray}
R \cdot R &=& L_{R}\, Q(g,R) ,
\label{pseudo}
\end{eqnarray}
where $L_{R}$ is some function on this set. We note that
${\mathcal{U}}_{S} \cup {\mathcal{U}}_{C} = {\mathcal{U}}_{R}$ (see, e.g. \cite{DGHHY}). 
We mention that \cite{DG1987} is the first paper, 
in which manifolds satisfying (\ref{pseudo}) were called pseudosymmetric manifolds. 
It is easy to check that (\ref{pseudo}) is equivalent on ${\mathcal{U}}_{R}$ to 
$(R - L_{R}\, G) \cdot (R  - L_{R}\, G) = 0$.
Evidently, every semisymmetric manifold is pseudosymmetric.
The converse statement is not true.
It seems that the Schwarzschild spacetime, the Kottler spacetime, the Reissner-Nordstr\o m spacetime, 
as well as some Friedmann-Lema{\^{\i}}tre-Robertson-Walker spacetimes are the ``oldest'' examples 
of non-semisymmetric pseudosymmetric warped product manifolds (see, e.g., \cite{{DHV2008}, {DVV1991}}).
Pseudosymmetric manifolds also are named {\sl Deszcz symmetric spaces} (see, e.g., \cite{Verstraelen}).
We also note that (\ref{pseudo}) implies
\begin{eqnarray}
R \cdot S &=& L_{R}\, Q(g,S) ,\ \ \ \ 
R \cdot C \ =\ L_{R}\, Q(g,C) .
\label{Weyl-pseudo-bis}
\end{eqnarray}
The conditions 
(\ref{pseudo}) and (\ref{Weyl-pseudo-bis})
are equivalent on the set ${\mathcal{U}}_{S} \cap {\mathcal{U}}_{C}$ 
of any warped product manifold $M_{1} \times_{F} M_{2}$, 
with $\dim\, M_{1} = \dim M_{2} = 2$ \cite{30}. 
A semi-Riemannian manifold $(M,g)$, $n \geq 3$, is called {\sl Ricci-pseudosymmetric} 
if the tensors $R \cdot S$ and $Q(g,S)$ are linearly dependent at every point of $M$
{\cite{{20}, {30}, {D 10}, {DGHSaw}, {P31}}}.
This is equivalent on ${\mathcal{U}}_{S}$ to 
\begin{eqnarray}
R \cdot S &=& L_{S}\, Q(g,S) ,
\label{Riccipseudo07}
\end{eqnarray}
where $L_{S}$ is some function on this set. 
As it was mentioned in Introduction, 
every warped product manifold $\overline{M} \times _{F} \widetilde{N}$
of an $1$-dimensional $(\overline{M}, \overline{g})$ manifold and
an $(n-1)$-dimensional Einstein semi-Riemannian manifold $(\widetilde{N}, \widetilde{g})$, $n \geq 3$, with a warping function $F$,
is a  quasi-Einstein manifold. Such warped products also are Ricci-pseudosymmetric manifolds,
see, e.g., {\cite[Section 1] {ChDGP}} and Example 4.1 of this paper.

A semi-Riemannian manifold $(M,g)$, $n \geq 4$, is said to be {\sl Weyl-pseudosymmetric} 
if the tensors $R \cdot C$ and $Q(g,C)$ are linearly dependent at every point of $M$
\cite{{30}, {D 10}, {DGHSaw}}. 
This is equivalent on ${\mathcal{U}}_{C}$ to 
\begin{eqnarray}
R \cdot C &=& L_{1}\, Q(g,C) ,
\label{Weyl-pseudo}
\end{eqnarray}
where $L_{1}$ is some function on this set. 
Using (\ref{eqn2.1}), we can check that 
on every Einstein manifold $(M,g)$, $n \geq 4$,
(\ref{Weyl-pseudo}) turns into $R \cdot R = L_{1}\, Q(g,R)$.
For a presentation of results on the problem 
of the equivalence of pseudosymmetry, Ricci-pseudosymmetry and Weyl-pseudosymmetry
we refer to {\cite[Section 4] {DGHSaw}}.  
A semi-Riemannian manifold $(M,g)$, $n \geq 4$, is said {\sl to have a pseudosymmetric Weyl conformal curvature tensor}
if the tensors $C \cdot C$ and $Q(g,C)$ are linearly dependent at every point of $M$ \cite{{37}, {30}, {D 10}}. 
This is equivalent on ${\mathcal U}_{C}$ to (\ref{4.3.012}), 
where $L_{C}$ is some function on this set. 
We note that (\ref{4.3.012}) is equivalent on ${\mathcal{U}}_{C}$ to 
$(C - L_{C}\, G) \cdot (C  - L_{C}\, G) = 0$.

As it was stated in \cite{30}, any warped product manifold $M_{1} \times_{F} M_{2}$, 
with $\dim\, M_{1} = \dim M_{2} = 2$, satisfies (\ref{4.3.012}). Thus in particular,
the Schwarzschild spacetime, the Kottler spacetime
and the Reissner-Nordstr\o m spacetime satisfy (\ref{4.3.012}).
Recently manifolds with pseudosymmetric Weyl tensor 
were investigated in \cite{{DGHHY}, {DeHoJJKunSh}}.
Warped product manifolds $\overline{M} \times _{F} \widetilde{N}$, of dimension $\geq 4$,
satisfying the condition (\ref{genpseudo01}) on 
${\mathcal U}_{C} \subset \overline{M} \times _{F} \widetilde{N}$,
where $L$ is some function on this set,
were studied in \cite{{DD 3}, {49}}.
In \cite{49} necessary and sufficient conditions for  
$\overline{M} \times _{F} \widetilde{N}$ to be a manifold satisfying (\ref{genpseudo01}) are given.
In particular, in that paper it was proved that 
any $4$-dimensional warped product manifold $\overline{M} \times _{F} \widetilde{N}$, 
with an $1$-dimensional base $(\overline{M},\overline{g})$, 
satisfies (\ref{genpseudo01}) {\cite[Theorem 4.1] {49}}.
For details about the pseudosymmetric, Ricci-pseudosymmetric and Weyl-pseudosymmetric manifolds
as well other conditions of this kind, named pseudosymmetry type curvature conditions, we refer to the papers:
\cite{{ChDGP}, {D 10}, {DGHSaw}, {DHV2008}, {HaVer}} 
and also references therein.

If $(M,g)$, $n \geq 4$, is an Einstein semi-Riemannian manifold 
then ${\mathcal{U}}_{R} = {\mathcal{U}}_{C}$ and (\ref{eqn2.1}) yields 
\begin{eqnarray}
C &=& R - \frac{\kappa }{(n-1)n}\, G .
\label{eqn2.1einstein}
\end{eqnarray}
\begin{thm}
If $(M,g)$, $n \geq 4$, is a pseudosymmetric Einstein semi-Riemannian manifold satisfying (\ref{pseudo}) 
on ${\mathcal{U}}_{R} \subset M$ then on this set we have
$R \cdot R - Q(S,R) = ( L_{R} - \frac{\kappa}{n} ) Q(g,C)$,
$C \cdot C = ( L_{R} - \frac{\kappa}{(n-1) n} ) Q(g,C)$ and
$C \cdot R + R \cdot C = Q(S,C) + ( 2 L_{R} - \frac{\kappa }{n-1} ) Q(g,C)$.
\end{thm}
{\bf{Proof.}} The second condition of our assertion was proved in {\cite[Theorem 3.1] {37}}. Further,
using (\ref{pseudo}) and (\ref{eqn2.1einstein}) we obtain $R \cdot C = L_{R}\, Q(g,C)$ and 
\begin{eqnarray*}
& & 
R \cdot R - Q(S,R) \ =\   
\left( L_{R} - \frac{\kappa }{n} \right)  Q(g,R - \frac{\kappa }{(n-1)n}\, G ) \ =\  \left( L_{R} - \frac{\kappa }{n} \right) Q(g,C)\, ,\\
& & 
C \cdot R + R \cdot C \ =\ ( R - \frac{\kappa }{(n-1)n}\, G) \cdot R + L_{R}\, Q(g,C) \\
& &=  \
R \cdot R -  \frac{\kappa }{(n-1)n}\, G \cdot R + L_{R}\, Q(g,C) 
\ = \ \left( L_{R} - \frac{\kappa }{(n-1)n} \right)  Q(g,R) + L_{R}\, Q(g,C)\\
& &
= \  
\left( 2 L_{R} - \frac{\kappa }{(n-1)n} \right) Q(g,C) \ =\ Q(S,C) + \left( 2 L_{R} - \frac{\kappa }{n-1} \right) Q(g,C) ,
\end{eqnarray*}
completing the proof.
\newline

In {\cite[Section 2] {Tashiro}} a class of $4$-dimensional Einstein Riemannian manifolds was defined and investigated. 
As it was stated in {\cite[Remark 5.1] {27}} those manifolds are pseudosymmetric. 
If a non-quasi-Einstein semi-Riemannian manifold $(M,g)$, $n \geq 4$,
satisfies on ${\mathcal{U}}_{S} \cap {\mathcal{U}}_{C} \subset M$
(\ref{pseudo}) and (\ref{genpseudo01})
or
(\ref{pseudo}) and (\ref{4.3.012}),
then (\ref{eq:h7a}) holds on this set ({\cite[Theorem 3.2 (ii)] {DY1994}}, {\cite[Lemma 4.1] {P43}}).
We also have the following converse statement.
\begin{thm} \cite{{DGHSaw}, {G5}}
If $(M,g)$, $n \geq 4$, is a semi-Riemannian manifold satisfying
(\ref{eq:h7a}) on ${\mathcal U}_{S} \cap {\mathcal U}_{C} \subset M$ 
then on this set we have
\begin{eqnarray*}
S^{2} &=& \alpha _{1}\, S + \alpha _{2} \, g ,
\ \
\alpha _{1} 
\ =\ 
\kappa + \frac{(n-2)\mu -1 }{\phi} ,
\ \
\alpha _{2}
\ =\
\frac{\mu \kappa + (n-1) \eta }{\phi } ,\\ 
R \cdot C &=& L_{R}\, Q(g,C),
\ \
L_{R}
\ =\
\frac{1}{\phi} \left( (n-2) (\mu ^{2} - \phi \eta) - \mu \right) ,\\
R \cdot R &=& L_{R}\, Q(g,R),
\ \ 
R \cdot S \ =\ L_{R}\, Q(g,S),\\
R \cdot R &=& Q(S,R) + L \, Q(g,C) ,
\ \ 
L
\ =\
L_{R} + \frac{\mu }{\phi }
\ =\
\frac{n-2}{\phi} (\mu ^{2} - \phi \eta),\\
C \cdot C &=& L_{C}\, Q(g,C) ,
\ \  
L_{C} 
\ =\ L_{R} + \frac{1}{n-2} (\frac{\kappa }{n-1} - \alpha _{1} ) ,\\
C \cdot R &=& L_{C}\, Q(g,R), \ \   
C \cdot S \ =\ L_{C}\, Q(g,S) ,
\end{eqnarray*}
\begin{eqnarray*}
R \cdot C - C \cdot R &=& \frac{1}{n-2}\, Q(S,R)
+ \left( \frac{(n-1) \mu - 1}{(n-2) \phi} + \frac{\kappa }{n-1} \right) Q(g,R)\\
& &
+ \frac{\mu ((n-1) \mu - 1) - (n-1) \phi \eta }{(n-2) \phi} \, Q(S,G) ,
\end{eqnarray*}
\begin{eqnarray*}
R \cdot C - C \cdot R &=& 
\left( \frac{1}{\phi} ( \mu - \frac{1}{n-2} ) + \frac{\kappa }{n-1} \right) Q(g,R)
+ \left( \frac{\mu}{\phi } ( \mu - \frac{1}{n-2}) - \eta \right) Q(S,G) .
\end{eqnarray*}
\end{thm}
\noindent
{\bf{Remark 3.1.}} 
Let the curvature tensor $R$ of a semi-Riemannian manifold $(M,g)$, $n \geq 4$, 
has the decomposition (\ref{eq:h7a}) on ${\mathcal U}_{S} \cap {\mathcal U}_{C} \subset M$. 
In {\cite[Lemma 3.2] {R102}} it was shown that 
the decomposition (\ref{eq:h7a}) is unique on this set.
\begin{prop} 
If $(M,g)$, $n \geq 4$, is a semi-Riemannian manifold satisfying
(\ref{eq:h7a}) on ${\mathcal U}_{S} \cap {\mathcal U}_{C} \subset M$
then (\ref{newRoter02}) and (\ref{newRoter01}) hold on this set. 
\end{prop} 
{\bf{Proof.}}
On ${\mathcal U}_{S} \cap {\mathcal U}_{C} \subset M$ we have
$C \cdot R = Q(S,C) + ( L_{R} - \frac{\kappa }{n-1}) Q(g,C)$  
{\cite[eq. (37)] {Kow02}}, 
where the function $L_{R}$ is defined by (\ref{newRoter02}) (see also Theorem 3.1). 
But this, together with $R \cdot C = L_{R}\, Q(g,C)$ 
and $L + L_{C} - \frac{1}{ (n-2) \phi } = 2 L_{R} - \frac{\kappa }{n-1}$
(see Theorem 3.1), completes the proof.
\begin{thm}
Let $(M,g)$, $n \geq 4$, be a semi-Riemannian manifold. 
(i)
The identity (\ref{identity01}) is satisfied on $M$.
(ii) 
If (\ref{genpseudo01}), 
with some function $L$,
is satisfied on ${\mathcal U}_{C} \subset M$ then (\ref{02identity01}) holds on this set.
(iii) 
If 
(\ref{genpseudo01}) and (\ref{4.3.012}),
with some functions $L$ and $L_{C}$,
are satisfied on ${\mathcal U}_{C} \subset M$ 
then (\ref{identity05}) holds on this set.
(iv) 
If $(M,g)$ is a non-Einstein and non-conformally flat semi-Riemannian manifold satisfying
on ${\mathcal U}_{S} \cap {\mathcal U}_{C} \subset M$ the conditions:
(\ref{quasi02}), and
(\ref{genpseudo01}) and
(\ref{4.3.012}), 
with some functions $L$ and $L_{C}$, 
then (\ref{identity05quasi}) holds on this set.
(v) 
The equation (\ref{Goedelidentity05}) is satisfied on the G\"{o}del spacetime.
\end{thm}
{\bf{Proof.}}
(i) 
We have (cf. {\cite[Section 1] {SawCM}})
\begin{eqnarray*}
(n-2)^{2} \, (C - R) \cdot (C - R) 
&=& ( g\wedge S - \frac{\kappa }{ n-1}\, G ) \cdot  ( g\wedge S - \frac{\kappa }{ n-1} \, G ) ,
\end{eqnarray*}
which yields
\begin{eqnarray*}
(n-2)^{2} \, (C \cdot C - R \cdot C - C \cdot R + R \cdot R ) 
&=&
( g\wedge S) \cdot ( g\wedge S)
- \frac{\kappa }{ n-1}\, G \cdot ( g\wedge S) .
\end{eqnarray*}
But this, by (\ref{identity21}), turns into (\ref{identity01}).
(ii) 
It is easy to see that (\ref{identity01}), 
by making use of (\ref{4.3.012}) and the identities  
(\ref{DS7new}) and (\ref{identity02}) turns into 
(\ref{02identity01}).
(iii) 
Relations (\ref{identity02}), (\ref{4.3.012}), 
(\ref{genpseudo01}) and (\ref{identity01}) yield 
\begin{eqnarray*}
& &
C \cdot R + R \cdot C  
\ =\ Q(S,R) + ( L +  L_{C}  )\, Q(g,C)
+ \frac{1}{(n-2)^{2}}\, Q( S^{2} - \frac{\kappa }{ n-1}\, S , G)\nonumber\\
&=& Q(S,C) + ( L +  L_{C} )\, Q(g,C) - \frac{1}{n-2}\, Q(g, \frac{1}{2}\, S \wedge S) 
+  \frac{1}{(n-2)^{2}}\, Q( S^{2} - \kappa \, S , G) ,
\end{eqnarray*}
which by (\ref{DS7}) turns into (\ref{identity05}).
(iii) It is easy to see that
the conditions 
(\ref{quasi02}), (\ref{quasi03}), (\ref{quasi03quasi03}),
(\ref{identity05}) and $Q(g,G) = 0$ lead to (\ref{identity05quasi}). 
(iv-v)
The Ricci tensor $S$ of the G\"{o}del spacetime $(M,g)$ satisfies 
$S = \kappa \, \omega \otimes \omega $, 
where $\omega$ is an $1$-form \cite{Godel}. From the last equation we get easily 
$S \wedge S = 0$ and $S^{2} = \kappa \, S$. 
It is also known that $R \cdot R = Q(S,R)$ and $C \cdot C = \frac{\kappa}{6}\, Q(g,C)$ hold on $M$
{\cite[Theorem 2] {DeHoJJKunSh}}. 
Now (\ref{identity05quasi}) yields (\ref{Goedelidentity05}).
Our theorem is thus proved.
\newline

\noindent
{\bf{Remark 3.2.}}
In {\cite[Section 4(v)] {DeHoJJKunSh}} it was shown that on the G\"{o}del spacetime
the tensors 
$R\cdot C$, $C\cdot R$, $Q(g,R)$, $Q(S,R)$, $Q(g,C)$ and $Q(S,C)$ are linearly dependent. 
\newline

We also have the following result.
\begin{prop}
{cf. \cite[Proposition 3.2, Theorem 3.3, Theorem 4.4] {DGHHY}}
If $(M,g)$, $n \geq 4$, is a semi-Riemannian manifold
satisfying on 
${\mathcal U}_{S} \cap {\mathcal U}_{C} \subset M$ the conditions
(\ref{genpseudo01}), (\ref{4.3.012}) and 
\begin{eqnarray}
R \cdot S  &=& Q(g,D) ,
\label{eqD}
\end{eqnarray}
where $D$ is a symmetric $(0,2)$-tensor, then (\ref{pseudo}) holds on this set.
Moreover, at every point of ${\mathcal U}_{S} \cap {\mathcal U}_{C}$
we have  
$\mathrm{rank}\, (S - \alpha _{1}\, g) = 1$ or
$\mathrm{rank}\, (S - \alpha _{1}\, g) \geq 2$ and  (\ref{eq:h7a}),
where 
$\alpha _{1} = \frac{1}{2} ( \frac{\kappa }{n-1} - L + L_{C})$.
\end{prop}
The last proposition, together with Proposition 3.3 and Theorem 3.4(iv), yields 
\begin{cor}
If $(M,g)$, $n \geq 4$, is a semi-Riemannian manifold
satisfying on 
${\mathcal U}_{S} \cap {\mathcal U}_{C} \subset M$ the conditions
(\ref{genpseudo01}), (\ref{4.3.012}) and (\ref{eqD}) then 
$\ C \cdot R  + R \cdot C \, =\, Q(S,C) + L_{2} \, Q(g,C)$ 
holds on ${\mathcal U}_{S} \cap {\mathcal U}_{C}$, 
where $L_{2}$ is some function on this set.
\end{cor}

Let $M$, $n = \dim M \geq 4$, be a connected hypersurface 
isometrically immersed in a semi-Riemannian space of constant curvature 
$N_{s}^{n+1}(c)$, with signature $(s,n+1-s)$, 
where $c = \frac{\widetilde{\kappa}}{n (n+1)}$ and 
$\widetilde{\kappa}$ is its scalar curvature. 
It is known that (\ref{genpseudo01}) holds on $M$. Precisely, 
\begin{eqnarray}
R \cdot R &=& Q(S,R) - \frac{(n-2) \widetilde{\kappa} }{n(n+1)}\, Q(g,C) 
\label{900ab}
\end{eqnarray}
on $M$ {\cite[Proposition 3.1] {DV}}. Now, as an immediate consequence of Theorem 3.3, we have 
\begin{thm}
Let $M$ is a hypersurface isometrically immersed in $N_{s}^{n+1}(c)$, $n \geq 4$.
Then
\begin{eqnarray}
C \cdot R  + R \cdot C 
&=& Q(S,C) 
- \frac{(n-2) \widetilde{\kappa} }{n(n+1)}\, Q(g,C) 
+ C \cdot C \nonumber\\
& & - \frac{1}{(n-2)^{2}}\, Q(g, \frac{n-2}{2}\, S \wedge S - \kappa \, g\wedge S + g \wedge S^{2} ) 
\label{02identity01hyper}
\end{eqnarray}
holds on $M$.
Moreover, if (\ref{4.3.012}) is satisfied on ${\mathcal U}_{S} \cap {\mathcal U}_{C} \subset M$ then on this set we have
\begin{eqnarray}
C \cdot R + R \cdot C &=& Q(S,C) + \left( L_{C} -  \frac{(n-2) \widetilde{\kappa} }{n(n+1)} \right) Q(g,C)\nonumber\\
& & - \frac{1}{(n-2)^{2}}\, Q( g, 
\frac{n-2}{2}\, S \wedge S - \kappa \, g \wedge S + g \wedge S^{2} ) .
\label{identity05hyper}
\end{eqnarray}
If $M$ is a quasi-Einstein hypersurface satisfying (\ref{quasi02}) and (\ref{genpseudo01}) 
on ${\mathcal U}_{S} \cap {\mathcal U}_{C}$ then on this set we have 
\begin{eqnarray}
C \cdot R + R \cdot C  &=& Q(S,C) + \left( L_{C} -  \frac{(n-2) \widetilde{\kappa} }{n(n+1)} \right) Q(g,C) .
\label{identity05quasihyper}
\end{eqnarray}
\end{thm}

It is known that every $2$-quasi-umbilical hypersurface 
in a semi-Riemannian space of constant curvature $N_{s}^{n+1}(c)$, $n \geq 4$, 
satisfies (\ref{4.3.012}) {\cite[Theorem 3.1] {P38}}. Now Theorem 3.4 yields
\begin{thm}
If $M$ is a $2$-quasi-umbilical hypersurface isometrically immersed in $N_{s}^{n+1}(c)$, $n \geq 4$,
then (\ref{identity05hyper}) holds on ${\mathcal U}_{S} \cap {\mathcal U}_{C} \subset M$.
\end{thm}

Let $M$ be an $n$-dimensional Chen ideal submanifold of codimension $m$ 
isometrically immersed in an Euclidean space
$\mathbb{E}^{n+m}$, $n \geq 4$, $m \geq 1$
\cite{{Ch01}, {Ch02}}.  
It is known that (\ref{genpseudo01}) 
and
(\ref{4.3.012}) 
hold on ${\mathcal U}_{C} \subset M$ 
({\cite[Theorem 1] {DP-TVZ-Bras}}, 
see also
{\cite[Section 6] {Chen04}} and {\cite[Section 3.1] {DP-TVZ}}).
Now Theorem 3.3(ii) yields
\begin{thm}
If $M$, $n \geq 4$, is a Chen ideal submanifold of codimension $m$, $m \geq 1$, 
isometrically immersed in an Euclidean space $\mathbb{E}^{n+m}$ 
then (\ref{identity05}) holds on this set.
\end{thm}

\noindent
{\bf{Remark 3.3.}} (i)  We refer to \cite{G6} for further results 
on quasi-Einstein hypersurfaces $M$ in $N_{s}^{n+1}(c)$, $n \geq 4$, satisfying (\ref{4.3.012}). 
(ii) We refer to \cite{DP-TVZ} for curvature properties of pseudosymmetry type of Chen ideal submanifolds
in an Euclidean space.  
(iii) From (\ref{900ab}) it follows that 
every Einstein hypersurface $M$ in $N_{s}^{n+1}(c)$, $n \geq 4$,
is a pseudosymmetric manifold satisfying 
(\ref{pseudo}) and  
$L_{R} = \frac{\kappa }{n} - \frac{(n-2) \widetilde{\kappa} }{n(n+1)}$ 
on ${\mathcal U}_{R} \subset M$ (cf. {\cite[Section 5.5] {D 10}}).
Now from Theorem 3.1 we have 
\begin{eqnarray*}
R \cdot C + C \cdot R &=& Q(S,C) +  \frac{n-2}{n} \left( \frac{\kappa }{n-1} - \frac{2 \widetilde{\kappa} }{ n+1} \right) Q(g,C)
\end{eqnarray*}
on ${\mathcal U}_{R}$.  
We refer to \cite{Magid} for examples of semisymmetric Einstein hypersurfaces 
in some semi-Riemannian spaces of constant curvature. 
(iv) Let $M$ be a hypersurface in $N_{s}^{n+1}(c)$, $n \geq 4$.
If at every point of ${\mathcal U}_{C} \subset M$ 
the tensor $H^{2}$, the square of the second fundamental tensor $H$ of $M$,
is a linear combination of $H$ and the metric tensor $g$ of $M$ then 
(\ref{4.3.012}) holds on  ${\mathcal U}_{C}$ (see, e.g., {\cite[Section 1] {P38}}).
Moreover, in view of Theorem 3.5(iii), (\ref{identity05}) 
is satisfied on ${\mathcal U}_{C}$. 

\section{Warped product manifolds}

Let now $(\overline{M},\overline{g})$ and $(\widetilde{N},\widetilde{g})$,
$\dim \overline{M} = p$, $\dim N = n-p$, $1 \leq p < n$, 
be semi-Riemannian manifolds
covered by systems of charts $\{ U;x^{a} \}$ and 
$\{ V;y^{\alpha } \} $,
respectively.
Let $F$ be a positive smooth function on $\overline{M}$. 
The warped product $\overline{M} \times _F N$ 
of $(\overline{M},\overline{g})$ and $(\widetilde{N}, \widetilde{g})$ 
is the product manifold $\overline{M} \times \widetilde{N}$ 
with the metric 
$g = \overline{g} \times _F \widetilde{g} =
{\pi}_1^{*} \overline{g} + (F \circ {\pi}_1)\, {\pi}_2^{*} \widetilde{g}$,
where 
${\pi}_1 : \overline{M} \times \widetilde{N} \longrightarrow \overline{M}$ and 
${\pi}_2 : \overline{M} \times \widetilde{N} \longrightarrow \widetilde{N}$ 
are the natural projections on $\overline{M}$ and $\widetilde{N}$, respectively {\cite{{BO 1}, {K 1}, {ON}}}.
Let 
$ \{ U \times V ; x^{1}, \ldots ,x^{p},x^{p+1} =  y^{1}, \ldots , x^{n} = y^{n-p} \} $ 
be a product chart for $\overline{M} \times \widetilde{N}$. The local
components $g_{ij}$ of the metric $g = \overline{g} \times _F \widetilde{g}$ with respect
to this chart are the following
$g_{ij} = \overline{g}_{ab}$ if $i = a$ and $j = b$,
$g_{ij} = F\, \widetilde{g}_{\alpha \beta }$ if $i = \alpha $ and $j = \beta $, and
$g_{ij} = 0$ otherwise,
where $a,b,c, d, f \in \{ 1, \ldots ,p \} $, 
$ \alpha , \beta , \gamma , \delta \in \{ p+1, \ldots ,n \} $
and $h, i, j, k, l, m, r, s \in \{ 1,2, \ldots ,n \} $. We will denote by bars
(resp., by tildes) tensors formed from $\overline{g}$ (resp., $\widetilde{g}$).
The local components 
\begin{eqnarray*}
\Gamma ^{h} _{ij} \ =\ \frac{1}{2}\, g^{hs} ( \partial_{i} g_{js} + \partial_{j} g_{is} - \partial_{s} g_{ij}),
\ \ 
\partial _j \ =\ \frac{\partial }{\partial x^{j}} ,
\end{eqnarray*} 
of the Levi-Civita connection $\nabla $
of $\overline{M} \times _F \widetilde{N}$ are the following 
\begin{eqnarray*}
\Gamma ^{a} _{bc}\, =\, \overline{\Gamma } ^{a} _{bc} , \ \
\Gamma ^{\alpha } _{\beta \gamma } 
\, =\, \widetilde{\Gamma } ^{\alpha } _{\beta \gamma } ,  \ \
\Gamma ^{a} _{\alpha \beta } 
\, =\, - \frac{1}{2} \bar{g} ^{ab} F_b \widetilde{g} _{\alpha \beta } , \ \
\Gamma ^{\alpha } _{a \beta } 
\, =\, \frac{1}{2F} F_a \delta ^{\alpha } _{\beta } , \ \
\Gamma ^{a} _{\alpha b}\ =\ \Gamma ^{\alpha } _{ab}\, =\, 0 , \ \
F_a\, =\, \frac{\partial F}{\partial x^{a}} 
\end{eqnarray*}
(see, e.g., \cite{{P54}, {K 2}}). The local components
\begin{eqnarray*}
R_{hijk}\ =\ g_{hs}R^{s}_{\, ijk}\ =\ g_{hs} (\partial _k \Gamma ^{s} _{ij} 
- \partial _j \Gamma ^{s} _{ik} + \Gamma ^{r} _{ij} \Gamma ^{s} _{rk}
- \Gamma ^{r} _{ik} \Gamma ^{s} _{rj} ) ,
\end{eqnarray*}
of the Riemann-Christoffel curvature tensor $R$
and the local components $S_{ij}$ of the Ricci tensor $S$
of the warped product $\overline{M} \times _F \widetilde{N}$ which may not vanish
identically are the following:
\begin{eqnarray}
R_{abcd} &=& \overline{R}_{abcd} ,\ \ \ 
R_{\alpha ab \delta} \ =\
- \frac{1}{2}\, T_{ab} \widetilde{g}_{\alpha \delta} ,\ \ \ 
R_{\alpha \beta \gamma \delta} \ =\
F \widetilde{R}_{\alpha \beta \gamma \beta} 
- \frac{1}{4}\, \Delta_1 F\, \widetilde{G}_{\alpha \beta \gamma \delta}\, ,
\label{L3}\\
S_{ab} &=& \overline{S}_{ab} - \frac{n-p}{2}\, \frac{1}{F}\, T_{ab} ,\ \ \  
S_{\alpha \beta } \ =\  
\tilde{S}_{\alpha \beta } 
- \frac{1}{2}\, ( \mathrm{tr}(T) + \frac{n-p-1}{2F} \Delta _1 F )\, \widetilde{g}_{\alpha \beta } ,
\label{AL2}
\end{eqnarray}
\begin{eqnarray}
T_{ab}&=& 
\overline{\nabla }_a F_b - \frac{1}{2F} F_a F_b ,\ \
\mathrm{tr}(T) \ =\ \overline{g}^{ab} T_{ab} , \ \ 
\Delta _{1} F \ =\ {\Delta}_{1\, \overline{g}} F\ =\ 
\overline{g}^{ab} F_a F_b ,
\label{AL3}
\end{eqnarray}
where $T$ is the $(0,2)$-tensor with the local components $T_{ab}$.
The scalar curvature $\kappa $ of
$\overline{M} \times _F \widetilde{N}$ satisfies the following
relation
\begin{eqnarray}
\kappa &=& \overline{\kappa } + \frac{1}{F}\, \widetilde{\kappa }
- \frac{n - p}{F}\, ( \mathrm{tr}(T) + \frac{n - p -1 }{4F} \Delta _1 F) .
\label{AL4}
\end{eqnarray}

Warped products play an important role in Riemannian geometry 
(see, e.g., \cite{Besse, BO 1, K 2, ON}) 
as well as in the general relativity theory 
(see, e.g., \cite{Blau, GrifPod, ON, Schmutzer}). 
Many well-known spacetimes of this theory, 
i.e. solutions of the Einstein field equations, are warped products, 
e.g. the Schwarzschild, Kottler, Reissner-Nordstr\o m, 
Reissner-Nordstr\o m-de Sitter, Vaidya, Vaidya-Kottler, 
Vaidya-Reissner-Nordstr\o m, Vaidya-Bonnor,
as well as Robert\-son-Walker spacetimes. We recall that 
a warped product $\overline{M} \times _F \widetilde{N}$ 
of an $1$-dimensional manifold 
$(\overline{M}, \overline{g})$, $\overline{g}_{11} = - 1$, 
and a $3$-dimensional Riemannian space of constant curvature $(\widetilde{N}, \widetilde{g})$, 
with a warping function $F$, is said to be a {\sl Robert\-son-Walker spacetime}
(see, e.g., \cite{GrifPod, ON, Schmutzer}). It is well-known 
that the Robertson-Walker spacetimes are conformally flat quasi-Einstein manifolds.
More generally, one also considers warped products 
$\overline{M} \times _F \widetilde{N}$ 
of $(\overline{M}, \overline{g})$, $\dim \, \overline{M} = 1$, 
$\overline{g}_{11} = - 1$, with a warping function $F$ 
and an $(n-1)$-dimensional Riemannian manifold $(\widetilde{N}, \widetilde{g})$, $n \geq 4$. 
Such warped products are called {\sl generalized Robertson-Walker spacetimes} 
\cite{ARS01, ARS02}. We mention that Einstein generalized Robertson-Walker 
spacetimes were classified in \cite{ARS02}.
Curvature conditions of pseudosymmetry
type on spacetimes have been considered among others in 
\cite{P119, ChDGP, 27,  37, P61, DHasKhamS,  DeHoJJKunSh, P54, DePlaScher, DVV1991, Kow01, Kow02}.
\newline
 
\noindent
{\bf{Example 4.1.}} 
The warped product manifold 
$\overline{M} \times _{F} \widetilde{N}$, 
of an $1$-dimensional manifold $(\overline{M},\overline{g})$, $\overline{g}_{11} = \pm 1$,
and an $(n-1)$-dimensional semi-Riemannian Einstein manifold 
$(\widetilde{N}, \widetilde{g})$, $n \geq 5$, 
which is not of constant curvature, with a warping function $F$,
satisfies on 
${\mathcal U}_{S} \cap {\mathcal U}_{C} \subset \overline{M} \times _{F} \widetilde{N}$:
\begin{eqnarray}
& &
R \cdot S \ =\ L_{S}\, Q(g,S), \ \ 
L_{S} \ = \ - \frac{\mathrm{tr} T}{2 F} ,\ \ 
\mathrm{rank} (S - \alpha \, g ) \ =\ 1,\ \   
\alpha \ =\ \frac{\kappa }{n-1} - L_{S} ,\nonumber\\
& &
(n-2)\, (R \cdot C - C \cdot R) \ =\ Q(S,R) - L_{S}\, Q(g,R),
\label{quasi10}
\end{eqnarray}
{\cite[Theorem 4.1] {ChDGP}}. Furthermore, 
using (\ref{quasi03}), (\ref{eqn2.1}), (\ref{DS7new}), (\ref{identity02}) and (\ref{quasi10}) we get
\begin{eqnarray*}
Q(g,R) &=& Q(g,C) - \frac{1}{n-2}\, Q(S,G) ,\\
Q(S,R) &=& Q(S,C) + \frac{1}{n-2} \left( \alpha - \frac{\kappa}{n-1} \right) Q(S,G) ,\\
(n-2)\, (R \cdot C - C \cdot R) &=& Q(S,C) - L_{S}\, Q(g,C) .
\end{eqnarray*}

Using Theorem 3.4(i)-(iii), 
{\cite[Theorem 4.1] {49}} and {\cite[Theorem 2] {43}} we obtain 
\begin{thm}
Let $\overline{M} \times _{F} \widetilde{N}$ be the warped product manifold 
of an $1$-dimensional manifold $(\overline{M},\overline{g})$, $\overline{g}_{11} = \pm 1$,
and a $3$-dimensional semi-Riemannian manifold $(\widetilde{N}, \widetilde{g})$.
If $(\widetilde{N}, \widetilde{g})$ is not a space of constant curvature then 
(\ref{genpseudo01}) and (\ref{02identity01}) 
hold 
on ${\mathcal U}_{C} \subset \overline{M} \times _{F} \widetilde{N}$.
Moreover, if $(\widetilde{N}, \widetilde{g})$ is a quasi-Einstein manifold then  
(\ref{4.3.012}) and (\ref{identity05}) 
hold 
on ${\mathcal U}_{S} \cap {\mathcal U}_{C} \subset \overline{M} \times _{F} \widetilde{N}$.
\end{thm}

The Ricci tensor of the following $3$-dimensional Riemannian manifolds $(\widetilde{N}, \widetilde{g})$: 
the Berger spheres, the Heisenberg group $Nil_{3}$, 
$\widetilde{PSL(2,{\mathbb{R}})}$ - the universal covering of the Lie group $PSL(2,{\mathbb{R}})$ 
and the Lie group $Sol_{3}$ {\cite[Section 3] {LVWang}}, 
a Riemannian manifold isometric to an open part of the Cartan hypersurface {\cite[Section 2] {DG90}}
and 
some three-spheres of Kaluza-Klein type {\cite[Theorem 2 $(ii)_{a}$] {CalvaPerr}} 
have exactly two distinct eigenvalues. 
Evidently, these manifolds are quasi-Einstein, and in a consequence,  
pseudosymmetric (see, e.g., {\cite[Theorem 1] {43}}). 
For further examples of $3$-dimensional quasi-Einstein manifolds
we refer to \cite{BDV2006} (Thurston geometries and warped product manifolds)
and \cite{Kowalski1993} (manifolds with constant Ricci principal curvatures).

Theorem 4.1 leads to the following result.
\begin{thm}
The conditions 
(\ref{genpseudo01}), (\ref{4.3.012}) and
(\ref{identity05})   
are satisfied 
on the warped product manifold
$\overline{M} \times _{F} \widetilde{N}$
of an $1$-dimensional manifold 
$(\overline{M},\overline{g})$, $\overline{g}_{11} = \pm 1$,
and the $3$-dimensional Riemannian manifold $(\widetilde{N}, \widetilde{g})$ such as: 
the Berger sphere, $Nil_{3}$, 
$\widetilde{PSL(2,{\mathbb{R}})}$, $Sol_{3}$,
a Riemannian manifold isometric to an open part of the Cartan hypersurface,
or some three-spheres of Kaluza-Klein type.  
\end{thm} 

Using {\cite[Theorem 4.2] {49}}, {\cite[Theorem 3.5] {25}} and {\cite[Theorem 3] {43}} we can prove 
\begin{thm}
If $\overline{M} \times _{F} \widetilde{N}$ is the warped product manifold 
of an $1$-dimensional manifold $(\overline{M},\overline{g})$, $\overline{g}_{11} = \pm 1$,
and an $(n-1)$-dimensional quasi-Einstein conformally flat
semi-Riemannian manifold $(\widetilde{N}, \widetilde{g})$, $n \geq 5$, 
then the conditions 
(\ref{genpseudo01}), (\ref{4.3.012}) and (\ref{identity05}) 
are satisfied 
on ${\mathcal U}_{C} \subset \overline{M} \times _{F} \widetilde{N}$.
\end{thm}

\vspace{2mm}

We mention that recently curvature conditions of pseudosymmetry type 
of four-dimensional Thurston geometries were investigated in \cite{Hasni}.

\section{Pseudosymmetric warped product manifolds}

In this section we present some results on pseudosymmetric warped product manifolds.
\begin{thm} \cite{38}
The Riemann-Christoffel curvature tensor $R$ of the warped product manifold $\overline{M} \times _{F} \widetilde{N}$,
with $\dim \overline{M} = p$, $\dim \widetilde{N} = n - p$, $1 \leq p \leq n-1$, $n \geq 3$,
satisfies (\ref{pseudo}), i.e. $R \cdot R = L_{R}\, Q(g,R)$, on some coordinate domain of a point 
$x \in {\mathcal U}_R \subset \overline{M} \times _{F} \widetilde{N}$ 
if and only if 
the following relations are satisfied on this set
\begin{eqnarray*}
(\overline{R} \cdot \overline{R} ) _{abcdef} 
&=& L_{R}\, Q(\overline{g} , \overline{R} ) _{abcdef},\\
%\label{pseudo01}\\
%\end{eqnarray}
%\begin{eqnarray}
2 H^{f} _d \overline{R} _{fabc} 
&=& \frac{1}{F} \, (T_{ac} H_{bd} - T_{ab} H_{cd} ),\\
%\label{pseudo02}\\
%\end{eqnarray}
%\begin{eqnarray}
H_{ad} \, (\widetilde{R} _{\delta \alpha \beta \gamma } - \frac{ \Delta _1 F }{4F} \, \widetilde{G} _{\delta \alpha \beta \gamma } )
&=& - \frac{1}{2} T^{f} _{d} H_{fa} \, \widetilde{G} _{\delta \alpha \beta \gamma } ,\\
%\label{pseudo03}\\
%\end{eqnarray}
%\begin{eqnarray}
(\widetilde{R} \cdot \widetilde{R} ) _{\alpha \beta \gamma \delta \lambda \mu } 
&=&
(F L_{R}  + \frac{ \Delta _1 F }{4F} )\, Q(\widetilde{g} , \widetilde{R} ) _{\alpha \beta \gamma \delta \lambda \mu } ,
%\label{pseudo03a}
\end{eqnarray*}
where $T_{ab}$ is defined by (\ref{AL3}) and
\begin{eqnarray}
H_{ab} &=& \frac{1}{2}\, T_{ab} + FL_{R}\, \overline{g} _{ab}\, .
\label{pseudo04}
\end{eqnarray}
\end{thm}
\begin{prop}
Let $\overline{M} \times _{F} \widetilde{N}$ 
be the warped product of semi-Riemannian manifolds
$( \overline{M},\overline{g})$ and $(\widetilde{N},\widetilde{g})$, 
$\dim \overline{M} = p$, $\dim \widetilde{N} = n - p$, $2 \leq p \leq n-1$, $n \geq 4$,
with the warping function $F$, 
and let $( \overline{M},\overline{g})$ and $(\widetilde{N},\widetilde{g})$ 
be a spaces of constant curvature, 
provided that $p \geq 3$ and $n - p \geq 3$, respectively.
(i) cf. {\cite[Corollary 2.1] {38}}
The warped product $\overline{M} \times _{F} \widetilde{N}$
satisfies $R \cdot R = L_{R}\, Q(g,R)$, i.e. (\ref{pseudo}), on some coordinate domain of a point 
$x \in {\mathcal U}_R \subset \overline{M} \times _{F} \widetilde{N}$ 
if and only if 
the following two relations are satisfied on this set 
\begin{eqnarray}
H_{ac} H_{bd} - H_{ab} H_{cd} 
&=& 
F \left( \frac{ \overline{\kappa } }{(p-1) p } - L_{R} \right) (\overline{g} _{ab} H_{cd} - \overline{g} _{ac} H_{bd} ) ,
\label{pseudo10}\\
%\end{eqnarray}
%\begin{eqnarray}
\widetilde{\kappa } \, H_{ad} 
&=&
(n-p)(n-p-1) 
\left(
( F L_{R} + \frac{ \Delta _1 F  }{4F} ) \, H_{ad}  - H^{2}_{ad} \right) ,
\label{pseudo11bb}
\end{eqnarray}
where 
$H^{2}_{ad} = \overline{g}^{ef} H_{ae}H_{df}$ and
$\overline{\kappa }$ and $\widetilde{\kappa }$ are the scalar curvatures of
$(\overline{M} ,\overline{g} )$ and $(\widetilde{N} ,\widetilde{g} )$, respectively.
Moreover, if $n-p \geq 2$ then (\ref{pseudo11bb}) is equivalent to
\begin{eqnarray}
H^{2}_{ad} 
&=& 
\left( \frac{F \overline{\kappa } }{(p-1)p } 
- \frac{ \widetilde{\kappa }}{(n-p-1) (n-p) } + \frac{ \Delta _1 F }{4F} \right) H_{ad} . 
\label{pseudo11b}
\end{eqnarray}
(ii)
If $H = \frac{ {{\rm tr}}\; H }{p} \, \overline{g}$ is satisfied on some coordinate domain $U$ of a point 
$x \in {\mathcal U}_R \subset \overline{M} \times _{F} \widetilde{N}$, for certain function $L_{R}$,  
then $T = \frac{ {{\rm tr}}\; T }{p} \, \overline{g}$ on $U$. 
Moreover, $H = 0$ and (\ref{pseudo}) hold on $U$,
provided that $L_{R} = - \frac{ {{\rm tr}}\; T }{2 p F}$. 
(iii)
Let $\overline{M} \times _{F} \widetilde{N}$ be the warped product 
satisfying $R \cdot R = L_{R}\, Q(g,R)$. If 
$H \neq \frac{ {{\rm tr}}\; H }{p} \, \overline{g}$ 
at a point $x \in {\mathcal U}_R \subset \overline{M} \times _{F} \widetilde{N}$ 
then on some neighbourhood $U \subset {\mathcal U}_{R}$ of $x$ we have 
\begin{eqnarray}
&(a)& L_{R}  \ =\  \frac{ \overline{\kappa } }{(p-1)p } ,\ \ \ \ (b)\ \ \mathrm{rank}\, (H) \ =\ 1 .
\label{pseudo14}
\end{eqnarray}
\end{prop}
{\bf{Proof.}} (i) This assertion is a consequence of Theorem 5.1 and the definition of
$H_{ab}$ (see (\ref{pseudo04})).
(ii) From our assumption, by (\ref{pseudo04}), it follows that $T$ is proportional to $g$. 
It is obvious that if we set $L_{R} = - \frac{ {{\rm tr}}\; T }{2 p F}$ then (\ref{pseudo04}) 
gives $H = 0$. Now (i) completes the proof of (ii).
(iii)
From (\ref{pseudo10}) we have
\begin{eqnarray*}
H_{cd} H_{ab} - H_{ac} H_{bd} 
&=& 
F \left( \frac{ \overline{\kappa } }{(p-1) p } - L_{R} \right) ( \overline{g} _{bd} H_{ac} - \overline{g} _{cd} H_{ab} ) ,
%\label{pseudo10}
\end{eqnarray*}
This, together with (\ref{pseudo10}), yields
\begin{eqnarray*}
\left( \frac{ \overline{\kappa } }{(p-1) p } - L_{R} \right) 
( \overline{g} _{ab} H_{cd} - \overline{g} _{ac} H_{bd}
+ \overline{g} _{bd} H_{ac} - \overline{g} _{cd} H_{ab} ) &=& 0 ,
%\label{pseudo10}
\end{eqnarray*}
which by contraction with $g^{ab}$ gives 
$( \frac{ \overline{\kappa } }{(p-1) p } - L_{R}) ( H_{cd} - \frac{{{\rm tr}}\; H }{p}\, \overline{g} _{cd} ) = 0$.
From this, by our assumption, we get immediately
(\ref{pseudo14})(a). Now (\ref{pseudo10}), by (\ref{pseudo14})(a), reduces to 
\begin{eqnarray}
H_{ac} H_{bd} - H_{ab} H_{cd} &=& 0 ,
\label{pseudo10dd}
\end{eqnarray}
which is equivalent to (\ref{pseudo14})(b).
Our proposition is thus proved.
\newline

%\noindent
Let $\overline{M} \times _{F} \widetilde{N}$ 
be the warped product of semi-Riemannian manifolds
$( \overline{M},\overline{g})$, $\dim \overline{M} = p$, and $(\widetilde{N},\widetilde{g})$, 
$\dim \widetilde{N} = n - p$, $2 \leq p \leq n-2$,
with the warping function $F$, 
and let $( \overline{M},\overline{g})$ and $(\widetilde{N},\widetilde{g})$ 
be spaces of constant curvature, 
provided that $p \geq 3$ and $n - p \geq 3$, respectively,
satisfying (\ref{pseudo})
on ${\mathcal U}_{R} \subset \overline{M} \times _{F} \widetilde{N}$. 
Moreover, let
$H \neq \frac{ {{\rm tr}}\; H }{p} \, \overline{g}$ at a point $x \in {\mathcal U}_R$.
We note that from (\ref{AL2}) it follows that $S \neq \frac{\kappa }{n}\, g$ at this point.
Further, in view of Proposition 5.2, 
(\ref{pseudo11b}), (\ref{pseudo14}) and (\ref{pseudo10dd})
hold on some neighbourhood $U \subset {\mathcal U}_{R}$ of $x$.
From (\ref{pseudo10dd}), by a suitable contraction, it follows that 
$H^{2} = {\mathrm{tr}} H \, H$ on $U$. The last equation and (\ref{pseudo11b}) yield
\begin{eqnarray}
\frac{\overline{\kappa} }{p }
+ \frac{ \widetilde{\kappa }}{(n-p-1) (n-p) F} 
+ \frac{ {{\rm tr}}\; T }{2 F} 
- \frac{ \Delta _1 F }{4F^{2}}  &=& 0 .
\label{pseudo21}
\end{eqnarray}
We note that if $p = 2$ then (\ref{WeylWeyl02}), (\ref{WeylWeyl06}), (\ref{WeylWeyl03}) and (\ref{pseudo21})
lead to $\rho _{0} = \rho = 0$, and $C = 0$.  
\newline

From the above presented considerations it follows 
\begin{thm}
Let $\overline{M} \times _{F} \widetilde{N}$ be the warped product
of a $2$-dimensional semi-Riemannian manifold $(\overline{M},\overline{g})$ and
an $(n-2)$-dimensional semi-Riemannian manifold $(\widetilde{N},\widetilde{g})$, $n \geq 4$, 
with the warping function $F$,  
and let $(\widetilde{N},\widetilde{g})$ be a space of constant curvature, 
provided that $n \geq 5$.
The manifold
$\overline{M} \times _{F} \widetilde{N}$ 
satisfies (\ref{pseudo}), i.e. $R \cdot R = L_{R} Q(g,R)$,
on some coordinate domain $U$ of a point 
$x \in {\mathcal U}_{S} \cap {\mathcal U}_{C} \subset \overline{M} \times _{F} \widetilde{N}$
if and only if  on $U$ we have
\begin{eqnarray*}
H &=& \frac{1}{2}\, T + F L_{R}\, \overline{g}
\ =\ \frac{1}{2}\, T + F ( - \frac{ {{\rm tr}}\; T }{4 F} ) \, \overline{g}
\ =\ \frac{1}{2}\, ( T - \frac{ {{\rm tr}}\; T }{2 }  \, \overline{g} )
\ =\ 0 , 
%\label{pseudo22}
\end{eqnarray*}
i.e. $T_{ab}$ is proportional to $\overline{g}_{ab}$ on $U$.
\end{thm}

At the end of this section we recall the following result of \cite{DeKow}.
\begin{thm}
cf. {\cite[Theorem 4.1] {DeKow}}
Let $\overline{M} \times _{F} \widetilde{N}$ be the warped product
of a $2$-dimensional semi-Riemannian manifold $(\overline{M},\overline{g})$ and
an $(n-2)$-dimensional semi-Riemannian manifold $(\widetilde{N},\widetilde{g})$, $n \geq 4$, 
with the warping function $F$,  
and let $(\widetilde{N},\widetilde{g})$ be a space of constant curvature, 
provided that $n \geq 5$. If the tensor $T_{ab}$, defined by (\ref{AL3}),
is proportional to $\overline{g}_{ab}$ at every point of 
${\mathcal U}_{S} \cap {\mathcal U}_{C} \subset \overline{M} \times _{F} \widetilde{N}$
then (\ref{eq:h7a}) holds on this set.
\end{thm}

\section{Warped product manifolds with $2$-dimensional base and Einsteinian fibre}

Let $\overline{M} \times _{F} \widetilde{N}$ be the warped product manifold 
of a $2$-dimensional manifold $(\overline{M},\overline{g})$
and an $(n-2)$-dimensional semi-Riemannian manifold $(\widetilde{N},\widetilde{g})$, $n \geq 4$, 
with a warping function $F$,
and let $(\widetilde{N},\widetilde{g})$ be an Einstein manifold,  
provided that $n \geq 5$. 
Now (\ref{AL2}) turns into
\begin{eqnarray}
S_{ad} &=& \frac{\overline{\kappa } }{2}\, g_{ab} - \frac{n-2}{2 F}\, T_{ab} ,\ \ \ 
S_{\alpha \beta } \ =\ \tau _{1} \, g_{\alpha \beta } ,
\label{RcciRicci01}\\
%\end{eqnarray}
%\begin{eqnarray}
\tau _{1} &=&  
\frac{\widetilde{\kappa} }{(n-2) F} - \frac{\mathrm{tr}(T) }{2F }  - (n-3)\, \frac{\Delta _1 F }{4 F^{2}}  .
\label{RcciRicci02}
\end{eqnarray}
From (\ref{RcciRicci01}) it follows that $T_{ab}$ is proportional to $\overline{g}_{ab}$ at a point of
${\mathcal U}_{S} \cap {\mathcal U}_{C} \subset \overline{M} \times _{F} \widetilde{N}$ if and only if
$S_{ab}$ is proportional to $\overline{g}_{ab}$ at this point. Furthermore, from (\ref{RcciRicci01}) also it follows that
(\ref{quasi0202weak}) is satisfied on
${\mathcal U}_{S} \subset \overline{M} \times _{F} \widetilde{N}$, i.e.
$\mathrm{rank}\, (S - \tau _{1} \, g ) \, \leq \, 2$ on this set.  
In addition, if 
$\mathrm{rank}\, (S - \tau _{1} \, g ) = 2$
then $\overline{M} \times _{F} \widetilde{N}$ is a $2$-quasi-Einstein manifold.
Thus we have
\begin{thm}
Let $\overline{M} \times _{F} \widetilde{N}$ 
be the warped product manifold
of a $2$-dimensional semi-Riemannian manifold $(\overline{M},\overline{g})$
and an $(n-2)$-dimensional semi-Riemannian manifold $(\widetilde{N},\widetilde{g})$, $n \geq 4$, 
with a warping function $F$, 
and let $(\widetilde{N},\widetilde{g})$ be an Einstein manifold, 
provided that $n \geq 5$. 
Then on ${\mathcal U}_{S} \subset \overline{M} \times _{F} \widetilde{N}$ we have
$\mathrm{rank}\, (S - \tau _{1} \, g ) \, \leq \, 2$,
where the function $\tau _{1}$ is defined by (\ref{RcciRicci02}).
Moreover, if $\mathrm{rank}\, (S - \tau _{1} \, g ) \, = \, 2$ on some open non-empty subset of ${\mathcal U}_{S}$
then $\overline{M} \times _{F} \widetilde{N}$ is a $2$-quasi-Einstein manifold.
\end{thm}

Let now $A$ be the $(0,2)$-tensor with the local components $A_{ij}$ defined by 
\begin{eqnarray}
A_{ij} &=& S_{ij} - \tau _{1}\, g_{ij} ,
\label{AA010bb}
\end{eqnarray}
where $\tau _{1}$ is the function defined by (\ref{RcciRicci02}). 
Using (\ref{RcciRicci01}) and (\ref{AA010bb}) we get
\begin{eqnarray}
A_{ad} &=&  S_{ad} - \tau _{1}\, g_{ad} ,\ \ \ 
A_{\alpha \beta }
\ =\ S_{\alpha \beta } - \tau _{1}\, g_{\alpha \beta } \ =\ 0 ,\ \ \
A_{a \alpha } \ =\ 0 .
\label{AA010}
\end{eqnarray}
From (\ref{AA010}) it follows immediately that 
$A_{ab}$ is proportional to $\overline{g}_{ab}$ 
if and only if $S_{ab}$ is proportional to $\overline{g}_{ab}$. 
Further, let $A^{2}$ be the $(0,2)$-tensor with the local components $A^{2}_{ij} = g^{rs}A_{ir}A_{js}$.
We have
\begin{eqnarray}
& &
A^{2}_{ij} \ =\ S^{2}_{ij} - 2 \tau _{1}\, S_{ij} + \tau _{1}^{2}\, g_{ij} ,\ \ 
A^{2}_{ad} \ =\ S^{2}_{ad} - 2 \tau _{1}\, S_{ad} + \tau _{1}^{2}\, g_{ad} ,\ \ 
A^{2}_{\alpha \beta } \ =\ 0,\ \ A^{2}_{a \alpha } \ =\ 0 ,
\label{AA010FF}
\end{eqnarray}
\begin{eqnarray}
& &
\mathrm{tr}(A) \ =\
g^{rs}A_{rs}
\ =\ \kappa - n  \tau _{1} ,\ \ 
\mathrm{tr}(A^{2}) \ =\
g^{rs}A^{2}_{rs} \ =\
\mathrm{tr}(S^{2}) - 2 \kappa \tau _{1} + n \tau _{1}^{2} ,\nonumber\\  
& &
\mathrm{tr}(A^{2}) - (\mathrm{tr}(A))^{2} \ =\ 
\mathrm{tr}(S^{2}) - \kappa ^{2} + (n-1) \tau _{1} ( 2 \kappa - n \tau _{1}) ,
\label{AA010DD}
\end{eqnarray}
and
\begin{eqnarray*}
& &
\mathrm{tr}(A^{2}) - (\mathrm{tr}(A))^{2} \ =\ g^{rs}A^{2}_{rs} - (g^{rs}A_{rs})^{2} \nonumber\\
&=& \overline{g}^{11}A^{2}_{11} + 2\, \overline{g}^{12}A^{2}_{12} + \overline{g}^{22} A^{2}_{22}
- ( \overline{g}^{11}A_{11} + 2\, \overline{g}^{12}A_{12} + \overline{g}^{22} A_{22} )^{2}\nonumber\\
&=&
\overline{g}^{11} g^{rs}A_{1r}A_{1s} + 2 \overline{g}^{12} g^{rs}A_{1r}A_{2s} + \overline{g}^{22} g^{rs}A_{2r}A_{2s} 
- ( \overline{g}^{11}A_{11} + 2\, \overline{g}^{12}A_{12} + \overline{g}^{22} A_{22} )^{2}\nonumber\\
&=&
\overline{g}^{11}\, ( \overline{g}^{11} (A_{11})^{2} + 2 \overline{g}^{12}A_{11}A_{12} + \overline{g}^{22}(A_{12})^{2} )
+
\overline{g}^{22}\, ( \overline{g}^{11} (A_{12})^{2} + 2 \overline{g}^{12}A_{12}A_{22} + \overline{g}^{22}(A_{22})^{2} )\nonumber\\
& &
+ 2\, \overline{g}^{12}\, ( \overline{g}^{11} A_{11}A_{12} + \overline{g}^{12}A_{11}A_{22} 
+ \overline{g}^{12}(A_{12})^{2} + \overline{g}^{22}A_{12} A_{22} )
- ( \overline{g}^{11}A_{11} + 2\, \overline{g}^{12}A_{12} + \overline{g}^{22} A_{22} )^{2}\nonumber\\
&=& - 2 \, ( \overline{g}^{11} \overline{g}^{22} -  ( \overline{g}^{12} )^{2} ) (A_{11} A_{22} - (A_{12})^{2} ) ,
\end{eqnarray*}
i.e. 
\begin{eqnarray}
\mathrm{tr}(A^{2}) - (\mathrm{tr}(A))^{2} &=& - 2 ( \mathrm{det} ( \overline{g} ))^{-1} \, (A_{11} A_{22} - (A_{12})^{2} ) .
\label{AA010abab}
\end{eqnarray}
From (\ref{AA010}) and (\ref{AA010abab}) it follows that 
at every point of $x \in {\mathcal U}_{S} \subset \overline{M} \times _{F} \widetilde{N}$ 
the conditions:
$\mathrm{rank}(A) = 2$ and  
$\mathrm{tr}(A^{2}) -  (\mathrm{tr}(A))^{2} \neq 0$ 
are equivalent.
Therefore on the set of at all points of ${\mathcal U}_{S}$ 
at which $\mathrm{rank}(A) = 2$
we can define the function $\tau _{2}$ by 
\begin{eqnarray}
\tau _{2}
&=& ( \mathrm{tr}(A^{2}) - (\mathrm{tr}(A))^{2})^{-1} .
\label{tau2}
\end{eqnarray}
We also note that in view of Lemma 2.2 we have
\begin{eqnarray}
A^{2}_{ad} &=& \mathrm{tr}(A)\, A_{ad} + \frac{1}{2} ( \mathrm{tr}(A^{2}) - (\mathrm{tr}(A))^{2})\, \overline{g}_{ad} ,\nonumber\\
Q(A,A^{2})_{abcd} 
&=& 
- \frac{1}{2} ( \mathrm{tr}(A^{2}) - (\mathrm{tr}(A))^{2}) \, Q(\overline{g},A)_{abcd} . 
\label{TT011}
\end{eqnarray}

We have
\begin{thm}
Let $\overline{M} \times _{F} \widetilde{N}$ 
be the warped product manifold
of a $2$-dimensional semi-Riemannian manifold $(\overline{M},\overline{g})$
and an $(n-2)$-dimensional semi-Riemannian manifold $(\widetilde{N},\widetilde{g})$, $n \geq 4$, 
with a warping function $F$, 
and let $(\widetilde{N},\widetilde{g})$ be an Einstein manifold, 
provided that $n \geq 5$. 
Moreover, let $V$ be the set of all points 
of ${\mathcal U}_{S} \cap {\mathcal U}_{C} \subset \overline{M} \times _{F} \widetilde{N}$ 
at which $S_{ad}$ is not proportional to $\overline{g}_{ad}$. Then on $V$ we have
\begin{eqnarray}
R \cdot S &=&
(\phi _{1} - 2 \tau_{1} \phi_{2} + \tau _{1}^{2} \phi_{3}) \, Q(g,S)
+ (\phi_{2} - \tau_{1} \phi_{3})\, Q(g, S^{2}) + \phi _{3}\, Q(S,S^{2}) ,
\label{genRps00}\\
\phi _{1} &=& \frac{2 \tau _{1} - \overline{\kappa} }{2 (n-2) } ,\ \ \ 
\phi _{2} = \frac{1}{n-2}, \ \ \
\phi _{3} = \frac{ \tau_{2} (2 \kappa - \overline{\kappa} - 2 (n-1) \tau_{1}) }{n-2} .
\label{genRps01}
\end{eqnarray}
The condition (\ref{Riccipseudo07})
holds on the set $({\mathcal U}_{S} \cap {\mathcal U}_{C}) \setminus V$.
\end{thm}
{\bf{Proof.}} 
Let $A$ be the $(0,2)$-tensor defined by (\ref{AA010}). Using now
(\ref{genformulas03}), (\ref{genformulas04}), (\ref{L3}), (\ref{AL2}), (\ref{RcciRicci01}) and (\ref{AA010}) we get 
\begin{eqnarray*}
- \frac{1}{2 F}\, T_{ad} 
&=& \frac{1}{n-2}\, ( A_{ad} + \frac{ 2 \tau_{1}  - \overline{\kappa}}{2}\, g_{ad} ) , \\
R_{abcd} 
&=& \frac{\overline{\kappa}}{2}\, G_{abcd} ,\ \ \ 
R_{a \alpha \beta d} \ =\ 
\frac{1}{n-2}\, ( A_{ad}  +  \frac{ 2 \tau _{1} - \overline{\kappa} }{2 (n-2) }\, g_{ad} )\, g_{\alpha \beta} ,
\end{eqnarray*}
\begin{eqnarray*}
(R \cdot A)_{abcd} &=& \frac{\overline{\kappa} }{2}\, Q(g,A)_{abcd} ,\ \ \
(R \cdot A)_{a \alpha \beta d} 
\ =\ \frac{1}{n-2}\, (A^{2}_{ad} +  \frac{ 2 \tau _{1} - \overline{\kappa} }{2 } \, A_{ad})\,  g_{\alpha \beta}, \\
(R \cdot A)_{\alpha \beta \gamma \delta } &=& 0 ,\ \ \
Q(g,A)_{\alpha \beta \gamma \delta } 
\ =\ Q(g,A^{2})_{\alpha \beta \gamma \delta }
\ =\ Q(A,A^{2})_{\alpha \beta \gamma \delta } \ =\ 0 .
\end{eqnarray*}
Using the above presented relations 
and (\ref{TT011}) we obtain on $V$ the following condition 
\begin{eqnarray}
R \cdot A &=& \phi _{1}\, Q(g,A) +  \phi _{2}\, Q(g,A^{2}) +  \phi _{3}\, Q(A,A^{2}) ,
\label{genRps03}
\end{eqnarray}
where $\phi _{1}$, $\phi _{2}$ and $\phi _{3}$ are defined on $V$ by (\ref{genRps01}).
Now (\ref{genRps03}), by (\ref{AA010}) and (\ref{AA010FF}), turns into (\ref{genRps00}). 
Finally, from (\ref{RcciRicci01}) and the fact that $S_{ab}$ is proportional to $g_{ab}$
on $({\mathcal U}_{S} \cap {\mathcal U}_{C}) \setminus V$ 
it follows that $T_{ab}$ also is proportional to $g_{ab}$ on this set.
Therefore, in view of {\cite[Corollary 4.1] {20}},
(\ref{Riccipseudo07}) holds on $({\mathcal U}_{S} \cap {\mathcal U}_{C}) \setminus V$. 
The last remark completes proof.

\section{Warped product manifolds with $2$-dimensional base and fibre of constant curvature}

We consider the warped product manifold 
$\overline{M} \times _{F} \widetilde{N}$ 
of a $2$-dimensional manifold $(\overline{M},\overline{g})$
and an $(n-2)$-dimensional semi-Riemannian manifold $(\widetilde{N},\widetilde{g})$, $n \geq 4$, 
with a warping function $F$,
and let $(\widetilde{N},\widetilde{g})$ be a space of constant curvature, 
provided that $n \geq 5$. 
Using Lemma 1.1, (\ref{L3})-(\ref{AL4}) 
and {\cite[eqs. (12)-(16)] {30}} 
we can check that the local components $C_{hijk}$ 
of the tensor Ricci tensor $S$ and the Weyl conformal curvature tensor $C$ 
of $\overline{M} \times _{F} \widetilde{N}$
are expressed by 
\begin{eqnarray}
C_{abcd} &=& \frac{(n-3) \rho_{0} }{n-1}\, G_{abcd} ,\ \ \
C_{\alpha bc \delta} \ =\ - \frac{(n-3) \rho_{0} }{(n-2) (n-1)}\, G_{\alpha bc \delta } ,\nonumber\\
C_{\alpha \beta \gamma \delta} &=&  \frac{2 \rho_{0} }{(n-2) (n-1)}\, G_{\alpha \beta \gamma \delta } ,\ \ \
C_{abc\delta }\ =\ C_{ab \gamma \delta } \ =\ C_{a \beta \gamma \delta } \ =\ 0 ,
\label{WeylWeyl01}
\end{eqnarray}
respectively, where
\begin{eqnarray}
\rho_{0} &=& 
\frac{ \overline{\kappa } }{2} + \frac{ \widetilde{\kappa } }{ (n-3)(n-2) F } 
+ \frac{\mathrm{tr}(T)}{2 F} - \frac{\Delta _1 F}{4 F^{2}} .
\label{WeylWeyl02}
\end{eqnarray}
We also have
\begin{eqnarray}
& &
F \tau_{1} + (n-3) \, \frac{\overline{ \kappa}}{2} F +  (n-2)\, \frac{\mathrm{tr}(T)}{2}  \nonumber\\ 
&=&
(n-3)\, \frac{\overline{ \kappa}}{2} F + \frac{ \widetilde{\kappa} }{n-2}
+ (n-3)\, \frac{\mathrm{tr}(T)}{2} - (n-3)\, \frac{ \Delta _{1} F }{4 F} \nonumber\\ 
&=& (n-3)\, F \left( \frac{\overline{ \kappa}}{2} + \frac{ \widetilde{\kappa} }{(n-3)(n-2)F} +  \frac{\mathrm{tr}(T)}{2 F}
- \frac{ \Delta _{1} F }{4 F^{2}} \right) 
\, =\, (n-3)\, F \rho_{0} \, =\, \frac{n-1}{2} \, F \rho ,
\label{newRcciRicci06new}
\end{eqnarray}
where
\begin{eqnarray}
\rho &=&  \frac{2 (n-3) \rho_{0} }{ n-1}  .
\label{WeylWeyl06}
\end{eqnarray}
Now the condition (\ref{WeylWeyl01}),
by (\ref{WeylWeyl06}), turns into
\begin{eqnarray}
C_{abcd} &=& \frac{\rho }{2}\, G_{abcd} ,\ \ 
C_{\alpha bc \delta} \ =\ - \frac{\rho}{2 (n-2)}\,  G_{\alpha bc \delta } ,\nonumber\\
C_{\alpha \beta \gamma \delta} &=& \frac{\rho }{(n-3) (n-2)}\,   G_{\alpha \beta \gamma \delta } ,\ \ 
C_{abc\delta } \ =\ C_{ab \gamma \delta } \ =\ C_{a \beta \gamma \delta } \ =\ 0 .
\label{WeylWeyl03}
\end{eqnarray}

\noindent
{\bf{Remark 7.1.}}
Let
$\overline{M} \times _{F} \widetilde{N}$ 
be the
warped product manifold 
of a $2$-dimensional manifold $(\overline{M},\overline{g})$
and an $(n-2)$-dimensional semi-Riemannian manifold $(\widetilde{N},\widetilde{g})$, $n \geq 4$, 
with a warping function $F$,
and let $(\widetilde{N},\widetilde{g})$ be a space of constant curvature, 
provided that $n \geq 5$. 
(i)
From (\ref{WeylWeyl01}) it follows immediately that 
the manifold $\overline{M} \times _{F} \widetilde{N}$ is conformally flat if and only if 
the function $\rho _{0}$, defined by (\ref{WeylWeyl02}), vanishes on $\overline{M}$. 
(ii) We refer to {\cite[Lemma 3.3, Lemma 4.1, Lemma 4.3] {DDV 2}}, {\cite[Example 5.4 (i), (ii)] {DD 3}} 
and {\cite[Sections 4 and 5] {K 3}} 
for examples of conformally flat warped product manifolds $\overline{M} \times _{F} \widetilde{N}$, with $\dim \overline{M} \geq 2$. 
(iii) Recently warped product spacetimes $\overline{M} \times _{F} \widetilde{N}$, with $\dim \overline{M} = \widetilde{N} = 2$,
satisfying curvature conditions of pseudosymmetry type were studied in \cite{DHasKhamS}. 

\begin{thm} 
Let $\overline{M} \times _{F} \widetilde{N}$ be the warped product manifold 
of a $2$-dimensional semi-Riemannian manifold $(\overline{M},\overline{g})$
and an $(n-2)$-dimensional semi-Riemannian manifold $(\widetilde{N},\widetilde{g})$, $n \geq 4$, 
with a warping function $F$,
and let $(\widetilde{N},\widetilde{g})$ be a space of constant curvature, 
provided that $n \geq 5$. 
\newline
(i) The following three conditions are satisfied on the set
${\mathcal U}_{C} \subset \overline{M} \times _{F} \widetilde{N}$:
\begin{eqnarray}
C \cdot C &=& - \frac{\rho }{2(n-2)}\, Q(g,C) ,
\label{WeylWeyl07}
\end{eqnarray}
where $\rho$ is defined by (\ref{WeylWeyl06}), 
(\ref{genpseudo01}) with the function $L$ be defined by 
\begin{eqnarray}
L &=&
- \frac{n-2}{ (n-1) \rho }   
\left(
\overline{\kappa} \left( \tau_{1} + \frac{ \mathrm{tr}(T) }{2F } \right)
+ \frac{n-3}{4 F^{2} }\, (\mathrm{tr}(T^{2}) - (\mathrm{tr}(T))^{2}) 
\right) ,
\label{newRicciRicci08new}
\end{eqnarray}
where $\tau _{1}$ is defined by (\ref{RcciRicci02}),
and (\ref{identity05}) with $L_{C} = - \frac{\rho }{2(n-2)}$ and $L$ defined by 
(\ref{newRicciRicci08new}).
\newline 
(ii) 
Let $V$ be the set of all points 
of ${\mathcal U}_{S} \cap {\mathcal U}_{C} \subset \overline{M} \times _{F} \widetilde{N}$ 
at which $S_{ad}$ is not proportional to $\overline{g}_{ad}$. Then on $V$ we have: 
\begin{eqnarray}
C &=& - \frac{ (n-1) \rho \tau _{2} }{(n-3)(n-2)}  
\left( \frac{ n-2 }{2}\,  S \wedge S - \kappa \, g \wedge S + g \wedge S^{2}
- \frac{ \mathrm{tr}(S^{2}) - \kappa^{2} }{n-1} \, G \right) ,
\label{mainTT032main}
\end{eqnarray}
\begin{eqnarray}
R \cdot C + C \cdot R 
&=& Q(S,C) + 
\left( L - \frac{\rho }{2 (n-2)} + \frac{n-3 }{ (n-2) (n-1) \rho \tau_{2} } \right) Q(g,C) ,
\label{identity06}
\end{eqnarray}
\begin{eqnarray}
C \cdot R &=& 
- \frac{ 1 }{(n-2)^{2}}\, Q(
( \frac{ \rho  }{2} + (n-1)\, \rho  \tau _{1}^{2} \tau _{2} )\, S
- (n-1)\, \rho \tau _{1} \tau _{2}\, S^{2} , G)\nonumber\\
& & - \frac{ (n-1)\, \rho \tau _{2} }{(n-2)^{2}}\, g \wedge  Q(S , S^{2})
- \frac{\rho }{2 (n-2)}\, Q(g,C) ,
\label{geneinst01}
\end{eqnarray}
\begin{eqnarray}
R \cdot C 
&=& Q(S,C) + 
\left( L + \frac{n-3 }{ (n-2) (n-1) \rho \tau_{2} } \right) Q(g,C)
+ \frac{ (n-1)\, \rho \tau _{2} }{(n-2)^{2}}\, g \wedge  Q(S , S^{2}) 
\nonumber\\
& & 
+ \frac{ 1 }{(n-2)^{2}}\, Q( ( \frac{ \rho  }{2} + (n-1)\, \rho  \tau _{1}^{2} \tau _{2} )\, S
- (n-1)\, \rho \tau _{1} \tau _{2}\, S^{2},G) .
\label{identity07}
\end{eqnarray}
On the set $({\mathcal U}_{S} \cap {\mathcal U}_{C}) \setminus V$ the Weyl tensor $C$
is expressed by a linear combination of the Kulkarni-Nomizu products 
$S \wedge S$, $g \wedge S$ and $g \wedge g$.
\end{thm}
{\bf{Proof.}} 
(i) Using (\ref{genformulas01}) and  (\ref{genformulas02})
we can verify that the local components 
$(C \cdot C)_{hijklm}$ and $Q(g,C)_{hijklm}$ 
of the tensors $C \cdot C$ and $Q(g,C)$
which may not vanish are those related to 
\begin{eqnarray}
\ \ \ \ \ \ (C \cdot C)_{\alpha abcd \beta } 
&=&
 - \frac{(n-1)\rho^{2}}{4 (n-2)^{2}}\, g_{\alpha \beta} G_{dabc},\ \ \
(C \cdot C)_{a \alpha \beta \gamma d \beta } 
\ =\
\frac{(n-1)\rho^{2}}{4 (n-2)^{2}(n-3)}\, g_{ad} G_{\delta \alpha \beta \gamma },
\label{WeylWeyl04}\\
%\end{eqnarray}
%\begin{eqnarray}
\ \ \ \ \ \ Q(g,C)_{\alpha abcd \beta } 
&=& 
\frac{(n-1)\rho }{2 (n-2)}\, g_{\alpha \beta} G_{dabc}, \ \ \ 
Q(g,C)_{a \alpha \beta \gamma d \delta } 
\ =\ 
- \frac{(n-1)\rho }{2 (n-2)(n-3)}\,  g_{ad} G_{\delta \alpha \beta \gamma }
\label{WeylWeyl05}
\end{eqnarray}
(cf. {\cite[eqs. (8)-(11)] {37}}). From (\ref{WeylWeyl04}) and (\ref{WeylWeyl05}) it follows that
(\ref{4.3.012})
holds on ${\mathcal U}_{C} \subset \overline{M} \times _F \widetilde{N}$,
where $L_{C} = - \frac{\rho }{2(n-2)}$ and $\rho$ is defined by (\ref{WeylWeyl06}).

We prove now that (\ref{genpseudo01}) is satisfied.
First of all, we recall that
necessary and sufficient conditions for
warped products of two semi-Riemannian spaces of constant curvature 
satisfying that condition are given in \cite{49}. 
In particular, when the base $(\overline{M},\overline{g})$ is a $2$-dimensional manifold, 
$(\widetilde{N},\widetilde{g})$ a space of constant curvature (when $n \geq 5$), 
then (\ref{genpseudo01}) holds on ${\mathcal U}_{C} \subset \overline{M} \times _{F} \widetilde{N}$ 
if and only if
\begin{eqnarray}
& &\left(
\left( 
\frac{\widetilde{\kappa} }{n-2} - \frac{1}{2} \left( \mathrm{tr}(T) + \frac{n-3}{2F} \Delta _1 F \right) \right)
\left( \frac{\overline{ \kappa} }{2} + \frac{L}{n-2} \right) 
+   \frac{n-3}{n-2} \, \frac{F L \overline{ \kappa} }{2} \right)
\overline{G}_{dabc} \nonumber\\ 
\ \ 
&=& 
\frac{n-3}{4 F}\, ( T_{ab}T_{cd} - T_{ac}T_{bd} ) 
- \left( \frac{\overline{ \kappa} }{4} + \frac{L}{2} \right)
( \overline{g}_{ab}T_{cd} + \overline{g}_{cd}T_{ab}
- \overline{g}_{ac}T_{bd} - \overline{g}_{bd}T_{ac} )
\label{newRcciRicci02new}
\end{eqnarray}
on ${\mathcal U}_{C}$ (cf., {\cite[Section 7, eq. (40)] {49}}).
Applying in (\ref{newRcciRicci02new}) the relation (\ref{RcciRicci02}) 
and the definitions of the tensors $g \wedge T$ and  $T \wedge T$  
we obtain
\begin{eqnarray*}
\left(
\left( \overline{\kappa}  + \frac{2 L}{n-2} \right) F \tau_{1} 
+   \frac{n-3}{n-2} \, F L \overline{ \kappa} \right)
\overline{G}_{dabc}  
&=& 
\frac{n-3}{2 F}\, \frac{1}{2}\, (T \wedge T)_{dabc} 
- \left( \frac{\overline{ \kappa} }{2} + L \right) (g \wedge T)_{dabc} .
%\label{newRcciRicci03new}
\end{eqnarray*}
Thie last equation, together with
\begin{eqnarray*}
(g \wedge T)_{1221} &=& \mathrm{tr}(T)\, G_{1221} \ =\  \mathrm{tr}(T)\, \det (\overline{g}) ,\\
\frac{1}{2}\, (T \wedge T) _{1221} &=&  T_{11}T_{22} - (T_{12})^{2} 
\ =\ 
- \frac{1}{2} \, \det (\overline{g})\, (\mathrm{tr}(T^{2}) - (\mathrm{tr}(T))^{2}) ,
\end{eqnarray*}
leads to 
\begin{eqnarray*}
& &\left(
\left( \overline{\kappa}  + \frac{2 L}{n-2} \right) F \tau_{1} 
+   \frac{n-3}{n-2} \, F L \overline{ \kappa} \right) \det (\overline{g})
\nonumber\\ 
&=& 
- \frac{n-3}{2 F}\, \frac{1}{2}\, (\mathrm{tr}(T^{2}) - (\mathrm{tr}(T))^{2}) \, \det (\overline{g})
- \left( \frac{\overline{ \kappa} }{2} + L \right) \mathrm{tr}(T)\, \det (\overline{g}) ,
\end{eqnarray*}
\begin{eqnarray*}
\left(
\left( \overline{\kappa}  + \frac{2 L}{n-2} \right) F \tau_{1} 
+   \frac{n-3}{n-2} \, F L \overline{ \kappa} \right) 
&=& 
- \frac{n-3}{2 F}\, \frac{1}{2}\, (\mathrm{tr}(T^{2}) - (\mathrm{tr}(T))^{2}) 
- \left( \frac{\overline{ \kappa} }{2} + L \right) \mathrm{tr}(T) ,
\end{eqnarray*}
\begin{eqnarray*}
\left( \frac{2 }{n-2} F \tau_{1} + \frac{n-3}{n-2} \, F \overline{ \kappa} +  \mathrm{tr}(T)  \right) L
&=& 
- F \overline{\kappa} \tau_{1} 
- \frac{n-3}{4 F}\, (\mathrm{tr}(T^{2}) - (\mathrm{tr}(T))^{2}) 
-  \frac{\overline{ \kappa} }{2} \mathrm{tr}(T) ,
%\label{newRcciRicci04new}
\end{eqnarray*}
\begin{eqnarray*}
& &
\left(  F \tau_{1} + (n-3) \, F \frac{\overline{ \kappa}}{2} +  (n-2) \frac{\mathrm{tr}(T)}{2}  \right) L\\
&=& 
- \frac{n-2}{2} \left(
F \overline{\kappa} \left( \tau_{1} + \frac{ \mathrm{tr}(T) }{2F } \right)
+ \frac{n-3}{4 F}\, (\mathrm{tr}(T^{2}) - (\mathrm{tr}(T))^{2}) 
 \right) .
\end{eqnarray*}
This, by making use of
(\ref{RcciRicci02}), (\ref{WeylWeyl02}) and (\ref{WeylWeyl06}), turns into
\begin{eqnarray*}
(n-1) \, \rho L &=&
- (n-2) \left(
\overline{\kappa} \left( \tau_{1} + \frac{ \mathrm{tr}(T) }{2F } \right)
+ \frac{n-3}{4 F^{2}}\, (\mathrm{tr}(T^{2}) - (\mathrm{tr}(T))^{2}) 
 \right) ,
%\label{newRcciRicci07new}
\end{eqnarray*}
which, together with (\ref{newRcciRicci06new}) and (\ref{WeylWeyl06}), yields (\ref{newRicciRicci08new}).
Now Theorem 3.3(ii) completes the proof of (i).
(ii) First we prove that the following relation is satisfied on $V$:
\begin{eqnarray}
\ \ \ \ C &=& \frac{\phi _{1}}{2}\, A \wedge A + \phi _{2}\, g \wedge A + \phi _{3}\, G + \phi _{4}\, g \wedge A^{2} \nonumber\\
&=&
- \frac{(n-1)\rho \tau _{2} }{(n-3) (n-2)} 
\left( \frac{ n-2 }{2}\, A \wedge A - \mathrm{tr}(A)\, g \wedge A +  g \wedge A^{2} - \frac{1}{ (n-1) \tau _{2} }\, G \right) ,
\label{TT016}
\end{eqnarray}
\begin{eqnarray}
\phi _{1} 
&=& - \frac{(n-1) \rho \tau _{2} }{n-3} , \ \ \ \ 
\phi _{2} \ =\ \frac{(n-1)\rho \tau _{2} \mathrm{tr}(A) }{(n-3) (n-2)} , \nonumber\\
\phi _{3} &=& \frac{\rho }{(n-3) (n-2)} , \ \ \ \
\phi _{4} \ =\ 
- \frac{(n-1)\rho \tau _{2} }{(n-3) (n-2)} , 
\label{TT017}
\end{eqnarray}
where the $(0,2)$-tensor $A$ is defined on $V$ by (\ref{AA010}).
Let $B$ be the $(0,4)$-tensor defined on $V$ by
\begin{eqnarray*}
B &=& C - \frac{\phi _{1}}{2}\, A \wedge A - \phi _{2}\, g \wedge A - \phi _{3}\, G - \phi _{4}\, g \wedge A^{2} ,
\end{eqnarray*}
where $\phi _{1}, \ldots , \phi _{4}$ are some functions on $V$. Evidently, $B$ is generalized curvature tensor.
Let $B_{hijk}$ be the local components of $B$. We have 
\begin{eqnarray*}
\ \ \ \
B_{hijk} &=& C_{hijk}
- \phi _{1}\, ( A_{hk}A_{ij} - A_{hj}A_{ik} )
- \phi _{2}\, ( g_{hk}A_{ij} + g_{ij}A_{hk} - g_{hj}A_{ik} - g_{ik}A_{hj} ) \nonumber\\
& &
- \phi _{3}\, ( g_{hk}g_{ij} - g_{hj}g_{ik} )
- \phi _{4}\, ( g_{hk}A^{2}_{ij} + g_{ij}A^{2}_{hk} - g_{hj}A^{2}_{ik} - g_{ik}A^{2}_{hj} ) .
%\label{nnAA010cdcd}
\end{eqnarray*}
It is clear that $B$ vanish at a point $x \in V$ if and only if  
\begin{eqnarray*}
C_{hijk} 
&=& 
\phi _{1}\, ( A_{hk}A_{ij} - A_{hj}A_{ik} )
+ \phi _{2}\, ( g_{hk}A_{ij} + g_{ij}A_{hk} - g_{hj}A_{ik} - g_{ik}A_{hj} ) \nonumber\\
& &
+ \phi _{3}\, ( g_{hk}g_{ij} - g_{hj}g_{ik} )
+ \phi _{4}\, ( g_{hk}A^{2}_{ij} + g_{ij}A^{2}_{hk} - g_{hj}A^{2}_{ik} - g_{ik}A^{2}_{hj} ) 
%\label{AA010cdcd}
\end{eqnarray*}
at $x$. We note that from
(\ref{AA010}) and (\ref{WeylWeyl01}) 
it follows immediately 
that the local components $B_{hijk}$ of the tensor $B$ 
which may not vanish identically are the following:
$B_{abcd}$, $B_{\alpha bc \delta}$ and $B_{\alpha \beta \gamma \delta}$.
Thus we see that $B = 0$ at $x$ if and only if 
\begin{eqnarray}
& &
\left( \frac{\rho }{2} - \phi _{3} \right) G_{abcd} \ =\ \phi _{1}\, ( A_{ad}A_{bc} - A_{ac}A_{bd} ) 
+ \phi _{2}\, ( g_{ad}A_{bc} + b_{bc}A_{ad} - g_{ac}A_{bd} - g_{bd}A_{ac} )\nonumber\\
& &
+ \phi _{4}\, ( g_{ad}A^{2}_{bc} + b_{bc}A^{2}_{ad} - g_{ac}A^{2}_{bd} - g_{bd}A^{2}_{ac} )
,\nonumber\\
& &
 - \left( \frac{\rho}{2 (n-2)} + \phi _{3} \right) g_{bc} g_{\alpha \delta} \ =\ 
\phi _{2} \, A_{bc} g_{\alpha \delta} 
+ \phi _{4} \, A^{2}_{bc} g_{\alpha \delta} 
,\nonumber\\
& &
\left( \frac{\rho }{(n-3) (n-2)} - \phi _{3} \right) G_{\alpha \beta \gamma \delta } \ =\  0 
\label{TT018}
\end{eqnarray}
at $x$.
Further, (\ref{TT018}),  by (\ref{TT011}), is equivalent to
\begin{eqnarray*}
& &
\phi _{1}\, ( A_{11}A_{22} - A_{12}A_{12} )
\ =\ \frac{(n-1) \rho }{2 (n-3)}\, G_{1221} 
,\nonumber\\
& &
- ( \phi_{2} + \mathrm{tr}(A) \phi_{4} ) \, A_{bc}
\ =\ 
\left( \frac{(n-1)\rho}{2 (n-3) (n-2)} 
+ \frac{1}{2} ( \mathrm{tr}(A^{2}) - (\mathrm{tr}(A))^{2}) \phi_{4}
\right) g_{bc}  ,\nonumber\\
& &
\phi _{3} \ =\ \frac{\rho }{(n-3) (n-2)} .
\end{eqnarray*}
But this, together with (\ref{AA010abab}), (\ref{tau2}) 
and the fact that $A_{ab}$ is not proportional to $g_{ab}$, 
leads immediately to (\ref{TT017}).

From 
(\ref{TT016}), by 
(\ref{AA010bb}),
(\ref{AA010FF}),
(\ref{AA010DD})
and
(\ref{TT017}), we get (\ref{mainTT032main}).
Now 
(\ref{identity05}), together with 
(\ref{newRicciRicci08new}) and (\ref{mainTT032main}),
yields (\ref{identity06}).
Using (\ref{eqn2.1}), (\ref{DS7}) and (\ref{WeylWeyl07}) we obtain
\begin{eqnarray*}
C \cdot R 
&=&
C \cdot \left( C + \frac{1}{n-2}\, g \wedge S - \frac{ \kappa }{(n-2)(n-1)}\, G \right) 
\ =\ 
C \cdot C  + \frac{1}{n-2}\, g \wedge (C \cdot S)\nonumber\\
&=& 
- \frac{ 1 }{(n-2)^{2}}\, 
g \wedge Q(g, ( \frac{ \rho  }{2} + (n-1)\, \rho  \tau _{1}^{2} \tau _{2} )\, S - (n-1)\, \rho \tau _{1} \tau _{2}\, S^{2})\nonumber\\
& & 
-  \frac{ (n-1)\, \rho \tau _{2} }{(n-2)^{2}} \, g \wedge Q(S , S^{2}) 
- \frac{\rho }{2 (n-2)}\, Q(g,C) 
\end{eqnarray*}
and in a consequence (\ref{geneinst01}). 
From 
(\ref{identity06})
and
(\ref{geneinst01})
it follows immediately (\ref{identity07}). 

From (\ref{RcciRicci01}) and the fact that $S_{ab}$ is proportional to $g_{ab}$
on $({\mathcal U}_{S} \cap {\mathcal U}_{C}) \setminus V$ 
it follows that $T_{ab}$ also is proportional to $g_{ab}$ on this set.
Therefore, in view of Theorem 5.4, (\ref{eq:h7a}) holds on $({\mathcal U}_{S} \cap {\mathcal U}_{C}) \setminus V$. 
Now using (\ref{eqn2.1})
and (\ref{eq:h7a}) we can express the tensor $C$ by a linear combination of the Kulkarni-Nomizu products 
$S \wedge S$, $g \wedge S$ and $g \wedge g$.
The last remark completes proof of (ii).
\newline
 
\noindent
{\bf{Remark 7.1.}} 
(i) Let the curvature tensor $R$ of a semi-Riemannian manifold $(M,g)$, $n \geq 4$, satisfies
\begin{eqnarray}
R &=& \frac{\phi _{1}}{2}\, S \wedge S + \phi _{2}\, g \wedge S + \phi _{3}\, G + \phi _{4}\, g \wedge S^{2}
\label{eqn0101.01special}
\end{eqnarray}
on ${\mathcal U}_{S} \cap {\mathcal U}_{C} \subset M$, 
where $\phi _{1}, \phi _{2}, \ldots , \phi _{4}$
are some functions on this set. 
Evidently, if (\ref{quasi03quasi03}) holds at a point of ${\mathcal U}_{S} \cap {\mathcal U}_{C}$  
then (\ref{eqn0101.01special}) reduces to (\ref{eq:h7a}) at this point.
We can prove that if the tensor $S^{3}$ is not a linear combination 
of $g$, $S$ and $S^{2}$ at a point ${\mathcal U}_{S} \cap {\mathcal U}_{C}$ then 
the decomposition (\ref{eqn0101.01special}) is unique at this point.
We also note that (\ref{eqn0101.01special}), by (\ref{eqn2.1}), yields 
\begin{eqnarray*}
C &=& \frac{\phi _{1}}{2}\, S \wedge S + \left( \phi _{2} - \frac{1}{n-2}\right) g \wedge S 
+ \left( \phi _{3} + \frac{\kappa }{(n-2)(n-1)} \right) G + \phi _{4}\, g \wedge S^{2} .
%\label{eqn0101.01special02}
\end{eqnarray*}
(ii) 
Warped product manifolds $\overline{M} \times _{F} \widetilde{N}$, $\dim \overline{M} = 1$,
satisfying (\ref{eqn0101.01special}) are investigated in  \cite{DGJP-TZ01}.
\newline

\noindent
{\bf{Example 7.1.}} (i) 
Let $\overline{M}_{1} = \{ (v,r) \in {\mathbb{R}}^{2}\, :\, r> 0 \}$, 
resp., $\overline{M}_{2} = \{ (u,r) \in {\mathbb{R}}^{2}\, :\, r> 0 \}$,  
be an open connected non-empty subset of ${\mathbb{R}}^{2}$ 
and let on $\overline{M}_{1}$, resp., $\overline{M}_{2}$, 
the metric tensor $\overline{g}_{1}$, resp., $\overline{g}_{2}$, be defined by
\begin{eqnarray*}
\overline{g}_{1 ab} dx^{a}dx^{b} \ =\ - f_{1}\, dv^{2} + 2\, dvdr ,\ \ 
\overline{g}_{2 ab} dx^{a}dx^{b} \ =\ - f_{2}\, du^{2} - 2\, dudr ,
\end{eqnarray*}
where 
$x^{1} = v$, $x^{2} = r$ and $f_{1} = f_{1}(v,r)$, 
resp., $x^{1} = u$, $x^{2} = r$ and $f_{2} = f_{2}(u,r)$,
is a smooth function on $\overline{M}_{1}$, resp., $\overline{M}_{2}$,
and $a,b = 1,2$.
We consider the warped product manifold
$\overline{M}_{i} \times _{F} \widetilde{N}$, $i = 1,2$, 
of $(\overline{M}_{i}, \overline{g}_{i})$, $i = 1,2$, 
and the $2$-dimensional standard unit sphere $(\widetilde{N},\widetilde{g})$
with the warping function $F = F(r) = r^{2}$.
(ii) (a) According to {\cite[Section 29.5.2] {Blau}}
(see also {\cite[Section 9.5] {GrifPod}}) 
the {\sl Vaidya metrics} 
form a simple class of timedependent generalizations 
of the Schwarzschild metric \cite{Vaidya01}.
They can be obtained from the Schwarzschild metric written in ingoing or outgoing
Eddington-Finkelstein coordinates by replacing the constant mass $m$ by a mass
function $m(v)$ or $m(u)$ depending on an advanced or retarded time coordinate.
The metrics of the warped products manifolds $\overline{M}_{i} \times _{F} \widetilde{N}$, $i = 1,2$, 
defined in (i), provided that  
$f_{1}(v,r) = 1 - \frac{2 m(v)}{r}$, resp.,
$f_{2}(u,r) = 1 - \frac{2 m(u)}{r}$, 
are the {\sl Vaidya metrics} (see, e.g., {\cite[eq. (29.15)] {Blau}} and {\cite[eq. (9.32)] {GrifPod}}).
(b) The metric of $\overline{M}_{1} \times _{F} \widetilde{N}$, resp.,
$\overline{M}_{2} \times _{F} \widetilde{N}$,
is called the {\sl generalized Vaidya ingoing metric}, 
resp., {\sl outgoing metric} (see, e.g. {\cite[eq. (39.16)] {Blau}}.
In particular, the metric of $\overline{M}_{1} \times _{F} \widetilde{N}$
with the function $f = f(v,r) = f_{1}(v,r)$ defined by
$f(v,r) = 1 - \frac{2 m(v)}{r} - \frac{\Lambda r^{2}}{3}, \Lambda \ =\ const.$,
$f(v,r) = 1 - \frac{2 m(v)}{r} - \frac{ q^{2}}{r^{2}}, q \ =\ const.$,
$f(v,r) = 1 - \frac{2 m(v)}{r} - \frac{ q^{2}(v)}{r^{2}}$,
respectively,
is named the {\sl Vaidya-Kottler}, the {\sl Vaidya-Reissner-Nordstr\o m} 
and the {\sl Vaidya-Bonnor ingoing metric}, respectively,
(see, e.g. {\cite[eqs. (39.18), (39.19), (39.20)] {Blau}}). 
(iii) (a)
For the manifold $\overline{M}_{1} \times _{F} \widetilde{N}$, 
with $f_{1}(v,r) = 1 - \frac{2 m(v)}{r}$, we have:
$S_{vv} = \frac{2 m'}{r^{2}}$, 
$m' = \frac{d m}{d v}$, $m = m(v)$, 
and $S_{hk} = 0$, if $h \neq v$ or $k \neq v$, 
$S^{2} = 0$, $\kappa = 0$, $S \cdot R =  \frac{2 m }{r^{3}} \, g\wedge S$, 	
$C \neq 0$, in particular $C_{v r r v } = - \frac{2 m }{ r^{3}}$. Moreover, 
\begin{eqnarray}
C \cdot C &=& 
R \cdot R - Q(S,R) \ =\
\frac{1}{2}\, ( R \cdot C + C \cdot R - Q(S,C) ) \ =\
- \frac{ m }{r^{3}} \, Q(g,C)\,  .	
\label{VaidyaVaidya01}
\end{eqnarray}
(b)
For the manifold $\overline{M}_{2} \times _{F} \widetilde{N}$, 
with $f_{2}(u,r) = 1 - \frac{2 m(u)}{r}$, we have:
$S_{uu} = - \frac{2 m'}{r^{2}}$, 
$m' = \frac{d m}{d u}$, $m = m(u)$,
and $S_{hk} = 0$, if $h \neq u$ or $k \neq u$, 
$S^{2} = 0$, $\kappa = 0$, $S \cdot R =  \frac{2 m }{r^{3}} \, g\wedge S$, 	
$C \neq 0$, in particular $C_{u r r u } = - \frac{2 m }{ r^{3}}$. 
Moreover, we also have
(\ref{VaidyaVaidya01}) (with $m = m(u))$. 
(iv) For the metric of the manifold $\overline{M}_{1} \times _{F} \widetilde{N}$, 
with the function $f = f(v,r) = f_{1}(v,r)$, we have
\begin{eqnarray*}
S_{v v} &=& \frac{1}{r^{2}} \left( 
 \frac{r^{2}}{2} f_{rr}'' + r f_{r}' - \frac{r}{f} f_{v}' \right) g_{v v} ,\ \
S_{v r} \ =\ -  \frac{1}{r^{2}} \left( \frac{r^{2}}{2} f_{rr}'' + r f_{r}' \right) g_{v r} , 
\nonumber\\
S_{\alpha \beta} 
&=& \tau_{1} \, g_{\alpha \beta} ,\ \
\tau_{1}  \ =\ 
\frac{1}{r^{2}} \left( -  r  f_{r}' - f + 1  \right) ,\ \
\kappa  \ =\ - \frac{2}{r^{2}} 
\left( \frac{r^{2}}{2} f_{rr}'' + 2 r f_{r}' + f - 1 \right) ,
\end{eqnarray*}
\begin{eqnarray*}
A_{v v} &=& 
S_{v v} - \tau_{1} \, g_{v v} \ =\ 
\frac{1}{ r^{2}} 
\left( \frac{r^{2}}{2} f_{rr}''  + 2 r f_{r}' + f  - 1 - \frac{r}{f} f_{v} ' \right) g_{v v} ,
\nonumber\\
A_{v r} &=& 
S_{v r} - \tau_{1} \, g_{v r} \ =\ \frac{1}{r^{2}} 
\left(
- \frac{r^{2}}{2} f_{rr}'' +  f - 1 \right) g_{v r} ,
\nonumber\\
A_{r r} &=& S_{r r} - \tau_{1} \, g_{r r} \ =\ 0 , \ \ \ \
A_{\alpha \beta} \ =\
S_{\alpha \beta} - \tau_{1}  \, g_{\alpha \beta}
\ =\ 0 ,
\end{eqnarray*}
where $f_{rr}'' = \frac{\partial^2 f}{\partial r^2}, f_{r}' = \frac{\partial f}{\partial  r}$
and $f_{v}' = \frac{\partial f}{\partial  v}$. 
We set 
$\tau_{3} = \frac{r^{2}}{2} f_{rr}'' -  r f_{r}' +  f - 1$. 
Now we can state that $\overline{M}_{1} \times _{F} \widetilde{N}$ is a conformally flat manifold
if and only if the function $\tau _{3}$ is a zero function.  
Furthermore, on $U_{C} \subset \overline{M}_{1} \times _{F} \widetilde{N}$ we have
(\ref{WeylWeyl03}) and (\ref{WeylWeyl07}), 
with $n = 4$ and $\rho = - \frac{2}{3} \tau _{3} r^{- 2}$, as well as (\ref{genpseudo01}) with 
$L = ( (  f - 1) f_{rr}'' - \frac{1}{2} \left( f_{r}' \right)^{2} ) \tau _{3}^{-1}$.
We also note that $\overline{M}_{1} \times _{F} \widetilde{N}$ is an Einstein manifold 
if and only if the function $f$ satisfies on  $\overline{M}_{1}$ the following system of differential equations
\begin{eqnarray*}
\frac{r^{2}}{2} f_{rr}'' - f + 1 \ =\ 0, &&
%\label{tau5a}\\
\frac{r^{2}}{2} f_{rr}'' + 2 r f_{r}'  
- \frac{ r}{f} f_{v}' +  f - 1 \ =\ 0 .
%\label{tau5b}
\end{eqnarray*}
It is easy to see that at every point of $\overline{M}_{1} \times _{F} \widetilde{N}$ we have
$\mathrm{rank} A = 2$ if and only if at every point of $\overline{M}_{1}$
we have
$\frac{r^{2}}{2} f_{rr}'' +  f - 1 \neq 0$.
Finally, $A_{ab}$ is proportional to $g_{ab}$ at a point of $\overline{M}_{1} \times _{F} \widetilde{N}$
if and only if at this point we have
$r f_{rr}''  + 2 f_{r}' - \frac{1}{f} f_{v}' = 0$.
\newline

\noindent
{\bf{Example 7.2.}} (i)
Let $\overline{M}_{1} = \{ (u,r) \in {\mathbb{R}}^{2}\, :\, 
r - 2 m > 0\ (\mbox{or}\ r - 2 m < 0) \}$ 
be an open connected non-empty subset of ${\mathbb{R}}^{2}$ and let on $\overline{M}_{1}$ 
the metric tensor $\overline{g}$ be defined by
\begin{eqnarray*}
\overline{g}_{uu} du^{2} + 2 \overline{g}_{ur} du dr + \overline{g}_{rr} dr^{2}
\ =\ - \exp(2 \beta) \, f\, du^{2} + 2 \exp(\beta)\, dudr ,
%\label{gsmm01}
\end{eqnarray*}
where 
$x^{1} = u$, $x^{2} = r$, $f = 1 - \frac{2 m}{r}$,
and $m = m(u,r)$ and $\beta = \beta(u,r)$ are some smooth functions on $\overline{M}$.
Further, let $\widetilde{g}$
be the standard metric on the $2$-dimensional unit sphere $\widetilde{N} = S^{2}(1)$. 
We denote by 
$g = \overline{g} \times _{F} \widetilde{g}$,
where $F = F(r) = r^{2}$, 
the warped product metric of $\overline{M} \times _{F} \widetilde{N}$.
The metric $g$ is said to be the {\sl general spherically symmetric metric} 
in advanced Eddington-Filkenstein coordinates,
see, e.g., {\cite[Section 4.1] {SenTorr}}. 
(ii) The local components of the Ricci tensor $S$ of $\overline{M} \times _{F} \widetilde{N}$
which may not vanish identically are the following
\begin{eqnarray*}
S_{u u} &=& 
\frac{1}{r^{2} ( r - 2 m )} 
( - 2 r \exp (- \beta ) m_{u}'  
 - ( - 3  r^{2} + 6  r  m ) \beta_{r}' m_{r}' 
 - ( - r^{2}  + 2 r m ) m_{rr}''
\nonumber\\
& & 
- (r^{3}  - 2 r^{2} m ) \exp (- \beta ) \beta_{ur}'' 
- (2 r^{2}  - 5 r m + 2 m^{2} ) \beta_{r}'
\nonumber\\
& &
- (r^{3}  - 4 r^{2} m + 4 r m^{2} ) 
(( \beta_{r}')^{2} + \beta_{rr}'') ) \, g_{uu} ,
\nonumber\\
S_{u r} &=& \frac{1}{r^{2}}  (  - r^{2}  \exp (- \beta ) \beta_{ur}'' 
+ 3 r m_{r}' \beta_{r}'
+ r m_{rr}''
+ ( - 2 r  + m ) \beta_{r}' 
+ ( - r^{2}  + 2 r m ) ( (\beta_{r}')^{2} + \beta_{rr}''))
g_{u r} ,
\nonumber\\
\ \
S_{r r} &=& \frac{2}{r} \beta_{r}' ,\ \
S_{\phi \phi } \ =\ \tau_{1}\, g_{\phi \phi},\ \ 
S_{\theta \theta } \ =\ \tau_{1}\, g_{\theta \theta}, \ \
\tau_{1}\ =\ 
\frac{1}{r^{2}} ( 2 m_{r}' - (r  - 2 m ) \beta_{r}' ) ,
%\label{gsmm02}
\end{eqnarray*}
where
$g_{\phi \phi} = r^{2}\, \widetilde{g}_{\phi \phi}$, $ \widetilde{g}_{\phi \phi} = 1$,
$g_{\theta \theta} = r^{2} \, \widetilde{g}_{\theta \theta}$, 
$\widetilde{g}_{\theta \theta} = \sin ^{2} \phi$
and
$m_{r}' = \frac{\partial m }{\partial  r} ,
m_{rr}'' = \frac{\partial^{2} m}{\partial r^{2}} ,
m_{u}'  = \frac{\partial m}{\partial  u}, 
\beta_{r}' = \frac{\partial \beta }{\partial  r} ,
\beta_{rr}'' = \frac{\partial ^{2} \beta }{\partial  r ^{2}},
\beta_{ur}'' = \frac{\partial^{2} \beta }{\partial  u\,\partial  r}$.
(iii) 
In the class of the general spherically symmetric metrics $g$ we also have non-Einstein metrics.
For instance, from the above formulas it follows immediately 
that the metrics $g$ with $S_{rr} \neq 0$, i.e. with 
$\beta_{r}' \neq 0$, are non-Einstein metrics. 
Moreover, for such metrics $S_{ab}$ are non-proportional to $g_{ab}$, $a,b = 1,2$. 
Some general spherically symmetric $g$ also are non-conformally flat metrics. 
Namely, the metrics $g$ satisfying 
\begin{eqnarray*}
& &
r^{3} ( \exp ( - \beta )  \beta_{ru}'' + (\beta _{r}')^{2} + \beta_{rr}'' ) 
- r^{2} (  m_{rr}'' + \beta_{r}' + 2 m  (\beta_{r}')^{2} + 2 m \beta_{rr}'' + 3  \beta_{r}' m_{r}' )\\
& &
- r (5m \beta_{r}' + 4 m_{r}') - 6 m \ =\ 0
\end{eqnarray*}
are non-conformally flat. This means that for some general spherically symmetric metrics $g$ 
the set $V$, defined in Theorem 7.1, is a non-empty subset
of ${\mathcal U}_{S} \cap {\mathcal U}_{C} \subset \overline{M} \times _{F} \widetilde{N}$.
\newline

\noindent
{\bf{Example 7.3.}} 
(i) Let $\overline{M} 
= \{ (t,r) \in {\mathbb{R}}^{2}\, :\, t > 0\ \mbox{and}\ r > 0 \}$ 
be an open connected non-empty subset of ${\mathbb{R}}^{2}$ and let on $\overline{M}$ 
the metric tensor $\overline{g}$ be defined by
$\overline{g}_{ab} dx^{a} dx^{b} = dt^{2} + R^{2}(t) dr^{2}$, $a,b = 1,2$,
where $x^{1} = t$, $x^{2} = r$, and $R = R(t)$ is a smooth positive (or negative) function on $\overline{M}$.
Let $\overline{M} \times _{F} \widetilde{N}$ be the warped product manifold 
of the manifold $(\overline{M},\overline{g})$
and the $2$-dimensional standard unit sphere $(\widetilde{N},\widetilde{g})$
with the warping function $F = F(t,r) = ( f(r) R(t))^{2}$, 
where $f = f(r)$ is a smooth positive (or negative) function on $\overline{M}$. 
We denote by $g = \overline{g} \times _{F} \widetilde{g}$ the metric of $\overline{M} \times _{F} \widetilde{N}$. 
We mention that the metric $g$ was considerd in {\cite[Section 4] {GHS01}}
(see also {\cite[Section 6] {GHS02ZW}}). 
(ii)
We set $\rho _{0} = (f  f_{rr}'' - ( f_{r}')^{2} + 1) (f R)^{-2}$, where 
$f_{r}' \ =\ \frac{d f }{d r}$ and $f_{rr}'' = \frac{d f_{r}'}{d r}$.
We can check that the Weyl conformal curvature tensor $C$ of $g$ is a zero tensor
if and only if $\rho _{0} = 0$ on $\overline{M}$. 
Further, we have $S_{12} = S_{21} = 0$,   
$S_{11} = \lambda _{1}\, g_{11}$, $S_{22} = \lambda _{2} \, g_{22}$,  
$S_{\alpha \beta} = \tau _{1}\, g_{\alpha \beta} = ( f R)^{2} \tau _{1}\, \widetilde{g}_{\alpha \beta}$, $\alpha , \beta = 3,4$, where
\begin{eqnarray*}
\lambda _{1} 
&=& - 3 R_{tt}'' R^{-1},\ \ R_{t}' \ =\ \frac{d R }{d t},\ \ R_{tt}'' \ =\ \frac{d R_{t}'}{d t}\, ,\ \
\lambda _{2} 
\ =\ - ( f R R_{tt}'' + 2 f ( R_{t}')^{2} + 2 f_{rr}'' ) f^{-1} R^{-2}\, ,\\
\tau _{1} 
&=& - (f^{2} R R_{tt}'' + 2 f^{2} ( R_{t}')^{2} + f f_{rr}'' + (f_{r}')^{2} - 1 ) ( f R)^{-2}
\ =\ \lambda _{2} + \rho _{0} 
\end{eqnarray*}
and $\widetilde{g}_{\alpha \beta}$ are the local components of the metric $\widetilde{g}$. 
(iii)
From (ii) it follows that $\lambda _{1} = \lambda _{2}$ if and only if
$R R_{rr}'' -  ( R_{r}')^{2} = c_{1}$ and $f_{rr}'' = c_{1} f$ and $c_{1} = const.$
on $\overline{M}$.
(iv)
If $\rho _{0}$ is non-zero
at a point of ${\mathcal U}_{C} \subset \overline{M} \times _{F} \widetilde{N}$
then in view of (ii) $S - \frac{\kappa}{4}\, g \neq 0$ at this point. 
Thus we have  ${\mathcal U}_{C} \subset {\mathcal U}_{S} \subset \overline{M} \times _{F} \widetilde{N}$.
Moreover, the following relations are satisfied on ${\mathcal U}_{C}$
\begin{eqnarray*}
& &
R\cdot R - Q(S,R) 
\ =\ - \frac{2}{3} \rho_{0}\, Q(g,C) \, ,
\ \ C \cdot C \ =\ - \frac{1}{6} \rho _{0}\, Q(g,C)\, ,\\
& &
R \cdot C + C \cdot R \ =\ Q(S,C) - \frac{1}{6} ( \kappa + 2 \rho_{0} )\, Q(g,C) \, .
\end{eqnarray*} 
(v) 
If $\lambda = \lambda _{1} = \lambda _{2}$ at a point of ${\mathcal U}_{C}$ then $S_{ab} = \lambda \, g_{ab}$, and by 
(\ref{AL2}), $T_{ab} = \frac{ {{\rm tr}}\; T }{2}\, g_{ab}$ at this point.
Let $V$ be the set of all points of ${\mathcal U}_{C}$ having this property.
From (iii) it follows that for some functions $f$ and $R$ the set $V$ is non-empty and in view of {\cite[Theorem 4.1] {DeKow}}
we can state that (\ref{eq:h7a}) holds on this set.
(vi)
If $\lambda _{1} \neq \lambda _{2}$ at a point of ${\mathcal U}_{C}$ then $S_{ab}$ is not proportional to $g_{ab}$ at this point.
Let $V$ be the set of all points of ${\mathcal U}_{C}$ having this property.
From (iii) it follows that for some functions $f$ and $R$ the set  
$V$ is non-empty and in view of Theorem 7.1(ii) we can state that
(\ref{mainTT032main})-(\ref{identity07}) hold on this set.

\vspace{4mm}

\noindent
\footnotesize{Ryszard Deszcz, Ma\l gorzata G\l ogowska and Jan Je\l owicki\\
Department of Mathematics\\
Wroc\l aw University of Environmental and Life Sciences\\ 
Grunwaldzka 53, 50-357 Wroc\l aw, Poland}\\
{\sf E-mail: Ryszard.Deszcz@up.wroc.pl}\ \ 
{\sf E-mail: Malgorzata.Glogowska@up.wroc.pl}\ \
{\sf E-mail: Jan.Je\l owicki@up.wroc.pl}
\newline

\noindent
\footnotesize{Georges Zafindratafa\\
Laboratoire de Math\'{e}matiques et Applications de Valenciennes\\
Universit\'{e} de Valenciennes et du Hainaut-Cambr\'{e}sis\\
Le Mont Houy, Lamath-ISTV2, 59313 Valenciennes Cedex 9, France}\\
{\sf E-mail: Georges.Zafindratafa@univ-valenciennes.fr}

\end{document}